\documentclass[12pt]{amsart}
\usepackage{amscd}
\usepackage{amsmath}

\newcommand{\C}{{\mathbb{C}}}        
\newcommand{\Q}{{\mathbb{Q}}}        
\newcommand{\Z}{{\mathbb{Z}}}        

\newcommand{\oddots}{{\mathinner{\mkern1mu\raise1pt\vbox{\kern7pt\hbox{.}}\mkern2mu\raise4pt\hbox{.}\mkern2mu\raise7pt\hbox{.}\mkern1mu}}}

\newcommand{\Fx}{{F^{\ast}}}

\newcommand{\supast}[1]{{{#1}^\ast}}

\newcommand{\qform}[1]{{\langle{#1}\rangle}}                   

\newcommand{\basemu}{\mathbf{\mu}}
\newcommand{\mmu}[1]{\basemu_#1}     

\DeclareMathOperator{\Spin}{Spin}           

\newcommand{\op}[1]{#1^\mathrm{op}}

\newcommand{\Int}{{\mathrm{Int}\,}}

\newcommand{\Gal}{{\mathrm{{\mathcal{G}}a\ell}\,}}

\DeclareMathOperator{\Id}{Id}

\newcommand{\End}[1]{{\mathrm{End}_{#1}}}
\newcommand{\EndF}{\End{F}}

\newcommand{\Aut}[1]{{{\mathrm {Aut}}_{#1}}}
\newcommand{\aut}[1]{\Aut{#1}}

\newcommand{\NrdA}{{\mathrm{Nrd}_A}}

\newcommand{\longto}{\longrightarrow}
\newcommand{\iso}{\xrightarrow{\sim}}
\newcommand{\injects}{\hookrightarrow}

\newcommand{\inv}[1]{{{#1}^{-1}}}

\newcommand{\Implies}{{\ \Longrightarrow\ }}
\newcommand{\If}{{\ \Longleftarrow\ }}

\newcommand{\QZt}{\mathbb{Q}/\mathbb{Z}(2)}
\newcommand{\oD}{{^1\hskip-.2em D_4}}

\newcommand{\g}{\mathfrak{g}}

\newcommand{\ozmat}{\tbtmat{1}{0}{0}{0}}
\newcommand{\sozmat}{\stbtmat{1}{0}{0}{0}}
\newcommand{\zomat}{\tbtmat{0}{0}{0}{1}}
\newcommand{\szomat}{\stbtmat{0}{0}{0}{1}}
\newcommand{\ubar}{\overline{u}}

\newcommand{\e}{\varepsilon}
\newcommand{\sym}{\mathcal{S}}
\newcommand{\s}{\sigma}

\newcommand{\ut}{{\underline{t}}}

\newcommand{\ul}{\underline{\lambda}}

\newcommand{\Ax}{\supast{A}}
\newcommand{\Dx}{\supast{\D}}

\newcommand{\Cd}{\C^d}
\newcommand{\Cdn}{(\Cd, \n)}

\newcommand{\bbrown}{(B, {-})}

\newcommand{\tbtmat}[4]{\left( \begin{array}{cc} #1&#2 \\ #3&#4 \end{array} \right) }
\newcommand{\stbtmat}[4]{\left( \begin{smallmatrix} #1&#2 \\ #3&#4 \end{smallmatrix} \right) }
\newcommand{\thbthmat}[9]{\left( \begin{array}{ccc} #1&#2&#3 \\ #4&#5&#6 \\ #7&#8&#9 \end{array} \right) }

\newcommand{\jordmat}[6]{\left( \begin{array}{ccc} #1&#6&\cdot \\ \cdot&#2&#4 \\ #5&\cdot&#3 \end{array} \right) }
\newcommand{\sjordmat}[6]{\left( \begin{smallmatrix} #1&#6&\cdot \\ \cdot&#2&#4 \\ #5&\cdot&#3 \end{smallmatrix} \right) }

\newcommand{\basjord}{\jordmat{\e_0}{\e_1}{\e_2}{a}{b}{c}}
\newcommand{\sbasjord}{\sjordmat{\e_0}{\e_1}{\e_2}{a}{b}{c}}

\newcommand{\sbasmat}{\stbtmat{\alpha}{j}{j'}{\beta}}
\newcommand{\basmat}{\tbtmat{\alpha}{j}{j'}{\beta}}

\newcommand{\trio}[1]{#1^{\times\negthinspace 3}}
\newcommand{\trion}[1]{(#1_0, #1_1, #1_2)}

\newcommand{\diag}{\mathrm{diag}\,}

\newcommand{\Inv}{\mathrm{Inv}\,}

\newcommand{\GInv}{\mathrm{GInv}\,}
\newcommand{\PInv}{\mathrm{PInv}\,}

\newcommand{\autp}{\Aut{}^{+}\,}
\renewcommand{\aut}{\Aut{}\,}
\newcommand{\GmF}{\mathbb{G}_{m,F}}

\newcommand{\Rel}{\mathrm{Rel}\,}

\newcommand{\SO}{{O^{+}}}

\newcommand{\OP}{\SO}
\newcommand{\GOP}{G\OP}

\newcommand{\B}{\mathcal{B}}

\renewcommand{\L}{\mathcal{L}}
\newcommand{\D}{\Delta}

\renewcommand{\C}{\mathfrak{C}}

\newcommand{\M}{\mathfrak{M}}

\newcommand{\h}{\mathfrak{H}}
\newcommand{\n}{\mathfrak{n}}

\newcommand{\Cn}{(\C,\n)}

\newtheorem{thm}{Theorem}[section]        
\newtheorem{lem}[thm]{Lemma}

\newtheorem{prop}[thm]{Proposition}

\newtheorem{dualitylem}[thm]{Duality Lemma}

\theoremstyle{definition}

\newtheorem{defn}[thm]{Definition}
\newtheorem{eg}[thm]{Example}

\theoremstyle{remark}

\newtheorem{rmk}[thm]{Remark}

\newenvironment{pf}{\noindent{\sl Proof:}}{\hfill$\qed$\medskip}

\numberwithin{equation}{section}

\makeatletter
\newenvironment{neqn}%
{\setcounter{equation}{\value{thm}}\begin{eqnarray}}%
{\end{eqnarray}{\stepcounter{thm}}\global\@ignoretrue}
\makeatother

\makeatletter
\newenvironment{ngather}%
{\setcounter{equation}{\value{thm}}}%
{{\stepcounter{thm}}\global\@ignoretrue}
\makeatother

\makeatletter
{\setcounter{equation}{\value{thm}}}%
{{\stepcounter{thm}}\global\@ignoretrue}
\makeatother

\makeatletter
\newenvironment{borel}[1]%
{\smallskip \noindent\refstepcounter{thm}{\bf \thethm.}{\bf{ #1.}}}%
{\smallskip \global\@ignoretrue}
\makeatother

\makeatletter
\newenvironment{borel*}%
{\smallskip \noindent\refstepcounter{thm}{\bf \thethm.}}%
{\smallskip \global\@ignoretrue}
\makeatother

\begin{document}

\title{Structurable algebras and groups of type $E_6$ and $E_7$}
\author{R.~Skip Garibaldi}
\date{\today}
\subjclass{17A30 (11E72 14L30 17A40 17B25 17C30 20G15)}

\begin{abstract}
It is well-known that every group of type $F_4$ is the automorphism
group of an exceptional Jordan algebra, and that up to isogeny all groups of 
type $^1E_6$ with trivial Tits algebras arise as the isometry groups
of norm forms  
of such Jordan algebras.  We describe a similar relationship between
groups of   
type $E_6$ and groups of type $E_7$ and use it to give explicit
descriptions of the  
homogeneous projective varieties associated to groups of type $E_7$
with trivial 
Tits algebras.
\end{abstract}

\maketitle

\nocite{Coop:FF1}
\nocite{Kantor}

It is well-known that over an arbitrary field $F$ (which for our
purposes we will assume has characteristic $\ne 2, 3$) every algebraic
group of type $F_4$ is obtained as the automorphism group of some
27-dimensional exceptional Jordan algebra and that some groups of type
$E_6$ can be obtained as automorphism groups of norm forms of such algebras.

In \cite{Brown:newalg}, R.B. Brown introduced a new kind of $F$-algebra,
which we will call a {\em Brown algebra}.  The automorphism groups of
Brown algebras provide a somewhat wider class of groups of type $E_6$,
specifically all of those with trivial Tits algebras.  Allison and
Faulkner \cite{AF:CD} showed that there is a
Freudenthal triple system (i.e., a quartic form and a skew-symmetric
bilinear form satisfying certain relations) determined up to
similarity by every Brown algebra.  The automorphism group of this
triple system is a simply 
connected group of type $E_7$, and we show that this provides a
construction of all simply connected groups of 
type $E_7$ with trivial Tits algebras (and more generally all
Freudenthal triple systems).  This is interesting because it 
allows one to relate properties of these algebraic groups over our
ground field $F$, which are generally hard to examine, with properties
of these algebras, which are relatively much easier to study.

Brown algebras are neither power-associative nor commutative, but they 
do belong to a wide class of algebras with involution known as central 
simple structurable algebras which were introduced in
\cite{A:structintro}, see \cite{A:survey} for a nice survey.
Other examples of such algebras are central
simple associative algebras with involution and Jordan algebras.
Brown algebras comprise the most poorly understood class of central
simple structurable algebras.  Although they are simple and 
so by definition have no 2-sided ideals which are stable under the
involution, they do have so-called ``inner ideals'', and we study them 
in Section \ref{idealsec}.  In particular, the largest inner ideal in
a Brown algebra is 12-dimensional and the largest singular ideal is
7-dimensional.

In Section \ref{flagsec}, we produce descriptions of the homogeneous
projective 
(a.k.a.~twisted flag) varieties associated to groups of type $E_7$
with trivial Tits algebras.  (These varieties are essentially the
spherical building associated to the group, see \cite[Ch.~V]{Brown} or
\cite[Ch.~5]{Ti:BN}.)
In another paper \cite{G:e7}, I define
objects called gifts (short for {\em g}eneralized {\em F}reudenthal
{\em t}riple system{\em s}) whose automorphism groups produce all groups of  
type $E_7$ over an arbitrary field up to isogeny.  The description of
the flag varieties here immediately gives a description of the
homogeneous projective varieties for arbitrary groups of type $E_7$ in
terms of 
gifts, which answers the question raised in \cite[p.~143]{MPW2}.

\subsection*{Notational conventions.}
All fields that we consider will have characteristic $\ne 2, 3$.  For
a field $F$, we write $F_s$ for its separable closure.

For $g$ an element in a group $G$, we write $\Int(g)$ for the
automorphism of $G$ given by $x \mapsto g x g^{-1}$.

For $X$ a variety over a field $F$ and $K$ any field extension of $F$,
we write $X(K)$ for the $K$-points of $X$.

When we say that an affine algebraic group (scheme) $G$ is {\em
  simple}, we mean that it is absolutely almost simple in the usual
sense (i.e., $G(F_s)$ has a finite center and no noncentral normal
subgroups).  For any simple algebraic group $G$ over a field $F$,
there is a unique minimal finite Galois field extension $L$ of $F$
such that $G$ is of inner type over $L$ (i.e., the absolute Galois
group of $L$ acts trivially on the Dynkin diagram of $G$).  We call
$L$ the {\em inner extension} for $G$.

We write $\GmF$ for the algebraic group whose $F$-points are $\Fx$ and
$\mmu{n}$ for the group of $n$th roots of unity. 

We will also follow the usual conventions for Galois cohomology and write
$H^i(F, G) := H^i(\Gal(F_s/F), G(F_s)$ for $G$ any algebraic group
over $F$, and similarly for the cocycles $Z^1(F, G)$.  For more
information about Galois cohomology, see \cite{SeLF} and \cite{SeCG}.

For $a, b \in \Fx$, we write $(a, b)_F$ for the (associative) {\em quaternion
$F$-algebra} generated by 
skew-commuting elements $i$ and $j$ such that $i^2 = a$ and $j^2 = b$,
please see \cite{Lam} or \cite[\S 14]{Draxl} for more information.

One oddity of the presentation should be pointed out to the reader.
We will be doing some explicit computations with Cayley algebras (see
\cite[Ch.~III, \S 4]{Schfr} or \cite[\S 33.C]{KMRT} for a definition), 
but not with their usual multiplication.  Instead, with juxtaposition
denoting the usual product and $\pi$ the standard involution, we will
use the multiplication $\star$ defined by $x \star y := \pi(x)
\pi(y)$.  This $\star$ multiplication is not even power-associative,
but it has some advantages when doing computations with exceptional
groups.  We will also make use of a particular basis $u_1, u_2,
\ldots, u_8$ of the {\em split Cayley algebra} $\Cd$, which has the
following multiplication table where each entry is $x \star y$ and
``$\cdot$'' replaces zero for clarity of reading:
\begin{neqn} \label{multtable}
\begin{tabular}{cc|cccc|cccc} 
   &     &\multicolumn{8}{c}{$y$} \\
   &     &  $u_1$&  $u_2$&  $u_3$&  $u_4$&  $u_5$&  $u_6$&  $u_7$&  $u_8$\\ \hline
   &$u_1$&$\cdot$&$\cdot$&$\cdot$& $-u_1$&$\cdot$& $-u_2$&  $u_3$& $-u_4$\\
   &$u_2$&$\cdot$&$\cdot$&  $u_1$&$\cdot$& $-u_2$&$\cdot$& $-u_5$& $-u_6$\\
   &$u_3$&$\cdot$& $-u_1$&$\cdot$&$\cdot$& $-u_3$& $-u_5$&$\cdot$&  $u_7$\\
$x$&$u_4$&$\cdot$& $-u_2$& $-u_3$&  $u_5$&$\cdot$&$\cdot$&$\cdot$& $-u_8$\\ \hline
   &$u_5$& $-u_1$&$\cdot$&$\cdot$&$\cdot$&  $u_4$& $-u_6$& $-u_7$&$\cdot$\\
   &$u_6$&  $u_2$&$\cdot$& $-u_4$& $-u_6$&$\cdot$&$\cdot$& $-u_8$&$\cdot$\\
   &$u_7$& $-u_3$& $-u_4$&$\cdot$& $-u_7$&$\cdot$&  $u_8$&$\cdot$&$\cdot$\\
   &$u_8$& $-u_5$&  $u_6$& $-u_7$&$\cdot$& $-u_8$&$\cdot$&$\cdot$&$\cdot$
\end{tabular}
\end{neqn}
For more discussion about this multiplication, please see \cite[\S 1]{G:iso}.

\section{Background on Albert algebras}

An {\em Albert algebra} over a field $F$ is a 27-dimensional central
simple exceptional Jordan algebra.  (Some of these adjectives are
redundant.)  Good introductions to Albert algebras may be found in
\cite{PR} or \cite[Ch.~IX]{Jac:J}, but we will recall what we need
here.

\begin{eg} \label{redalbeg}
Let $\C$ be a Cayley $F$-algebra and let $\gamma \in GL_3(F)$ be a
diagonal matrix.  Let $\ast$ denote the conjugate transpose on
$M_3(\C)$.  We write $\h_3(\C, \gamma)$ for the subspace of
$M_3(\C)$ fixed by $\Int(\gamma) \circ \ast$ and endowed with a
symmetrized product $\cdot$ given by
\[
a \cdot b := \frac12 (ab + ba),
\]
where juxtaposition denotes the usual product in $M_3(\C)$.  Then
$\h_3(\C, \gamma)$ is an Albert $F$-algebra.
\end{eg}

An Albert algebra is called {\em split} if it is isomorphic to
$\h_3(\C,1)$ for $\C$ the split Cayley algebra.  It is called reduced
if it is isomorphic to one
as in the preceding example.  We will want to do some explicit
computations in reduced Albert aglebras $\h_3(\C, \gamma)$ for $\gamma 
= \diag(\gamma_0, \gamma_1, \gamma_2)$.   For simplicity of notation
we will write 
\begin{neqn} \label{rednot}
\basjord 
\hbox{\ instead of\ }
\thbthmat{\e_0}{c}{\gamma_0^{-1} \gamma_2 \overline{b}}{\gamma_1^{-1} \gamma_0
  \overline{c}}{\e_1}{a}{b}{\gamma_2^{-1} \gamma_1 \overline{a}}{\e_2},
\end{neqn}
since the entries we have replaced with a $\cdot$ are forced by the
fact that elements of $\h_3(\C, \gamma)$ are fixed by the involution.

Every Albert $F$-algebra $J$ is endowed with a cubic norm map $N \!: J
\longto F$ and a linear trace map $T \!: J \longto F$.  We also use
$T$ to denote the map $T \!: J \times J \longto F$ given by
\[
T(x, y) := T(xy).
\]
This is a symmetric bilinear form on $J$, and it is nondegenerate 
\cite[p.~240, Thm.~5]{Jac:J}.  For any $f \in \EndF(J)$, we denote the
adjoint of 
$f$ with respect to $T$ by $f^\ast$, so $T(fx,y) = T(x, f^\ast y)$ for
all $x, y \in J$.  For simplicity, if $f$ is invertible we write
$f^\dagger$ for 
$(f^{-1})^\ast = (f^\ast)^{-1}$.

Now if one extends scalars to $F(t)$, expands $N(x + t y)$ and
considers the coefficient of $t$, then for a fixed $x \in J$ this
provides a linear map $f_x \!: J \longto F$ given by substituting in
for $y$.  Since $T$ is nondegenerate, there is an element $x^{\#} \in
J$ such that $T(x^{\#}, y) = f_x(y)$ for all $y \in J$.  Then $\#$
provides a quadratic map 
$\#\!: J \longto J$, cf.~\cite[pp.~495, 496]{McC:FST}.  We define a
linearization of $\#$ called the {\em Freudenthal cross product} by 
\[
x \times y := (x+y)^{\#} - x^{\#} - y^{\#},
\]
as in \cite[p.~496]{McC:FST}.  Note that $x \times x = 2x^{\#}$, which
differs by a factor of 2 from the definition of $\times$ given in
\cite{Jac:J}.

\begin{defn}
Let $J$ be an Albert algebra over $F$.  We call an $F$-vector space
map $\varphi \!: J \longto J$ a {\em norm similarity} if there is
some $\lambda \in \Fx$ such that $N(\varphi(j)) = \lambda N(j)$ for
every $j \in J$.  In that case, we call $\lambda$ the {\em multiplier}
of $\varphi$.  If $\varphi$ is a norm similarity with multiplier 1, we
call $\varphi$ a {\em norm isometry}.
\end{defn}

\subsection*{Algebraic groups}

For $J$ an Albert $F$-algebra, we define $\Inv(J)$ to be the algebraic
group whose $F$-points are the norm isometries of $J$.  (This is
Freudenthal's notation from \cite[\S 1]{Frd:E7.1}.)  This is a
simple simply connected algebraic group of type
$^1E_6$ over $F$.

Associated to any semisimple algebraic group are its {\em Tits
  algebras}, which are the endomorphism rings of irreducible
  representations, see \cite[\S 27]{KMRT} or \cite{Ti:R}.  In general,
  they are central simple algebras over finite separable extensions of
  $F$.  We say that a group has {\em trivial Tits algebras} if all of
  them are split.

There is a strong connection between Albert $F$-algebras and groups of
type $^1E_6$ over $F$ with trivial Tits algebras.  We summarize this
in the following theorem.  In \ref{isosim} and \ref{e7thm} we will
show that there is a 
similar connection between Brown $F$-algebras (defined in \ref{browndef}) and
groups of type $E_7$ over $F$ with trivial Tits algebras.

\begin{thm}
\begin{enumerate}
\item Every simple simply connected group of type
  $^1E_6$ over $F$ with trivial Tits algebras is isomorphic to
  $\Inv(J)$ for some Albert $F$-algebra $J$.

\item For $J_1$, $J_2$ Albert $F$-algebras, the following are
  equivalent:
\begin{enumerate}
\item $\Inv(J_1) \cong \Inv(J_2)$
\item $J_1$ and $J_2$ have similar norm forms
\item $J_1 \sim J_2$ (i.e., $J_1$ is isotopic to $J_2$, see below).
\end{enumerate}
\end{enumerate}
\end{thm}

In the statement of the preceding theorem, we used the notion of
isotopy of Jordan algebras which provides an equivalence relation for
such algebras which is weaker than isomorphism.  Specifically, for $u \in J$, 
we define a new Jordan algebra $J^{\qform{u}}$ which has the same underlying 
vector space as $J$ and whose multiplication $\cdot_u$ is given by
\begin{neqn} \label{jordisodef}
x \cdot_u y := \{ a, u, b \} = (a \cdot u) \cdot b + (b \cdot u) \cdot a - (a \cdot b) \cdot u,
\end{neqn}
where $\cdot$ denotes the usual multiplication in $J$.  We say that another 
Jordan algebra $J'$ is {\em isotopic} 
to $J$ (written $J' \sim J$) if $J'$ is isomorphic to $J^{\qform{u}}$
for some $u \in J$.  

\medskip
\begin{pf}
(1) follows from \cite[6.4.2]{Ti:R}.  That (2b) implies (2a) is clear,
and the converse is \cite[p.~38, Thm.~7]{Jac:ex}.  Finally, (2a) is
equivalent to (2c) by \cite[p.~55, Thm.~10]{Jac:ex}.
\end{pf}

\subsection*{Useful lemmas}

A very useful fact for us is that if $J$ is a reduced Albert
$F$-algebra, then there is a norm similarity of $J$ with multiplier
$\lambda$ for every $\lambda \in \Fx$.  Such a similarity is given by
$\psi$ for 
\begin{neqn} \label{jordmult}
\psi\basjord = 
\jordmat{\lambda \e_0}{\lambda \e_1}{\lambda^{-1} \e_2}{a}{b}{\lambda c}.
\end{neqn}

\begin{lem} \label{jordsim}
Suppose that $\varphi$ is a norm similarity of an Albert $F$-algebra
$J$ with multiplier $\lambda$.  Then $\varphi^\dagger$ is a norm
similarity for $J$ with multiplier $1/\lambda$, 
\[
\varphi(j_1) \times \varphi(j_2) = \lambda \varphi^\dagger (j_1
\times j_2), \hbox{\ and\ } 
\varphi^\dagger(j_1) \times \varphi^\dagger(j_2) =
\frac1{\lambda}\varphi(j_1 \times j_2).
\]
\end{lem}

\begin{pf}
Since these formulas hold if and only if they hold over
a field extension of $F$, we may assume that $F$ is algebraically
closed.  Let $\ell \in \Fx$ be such that $\ell^3 = \lambda$.  Then
$\varphi_1 := \frac1{\ell} \varphi$ is a norm isometry of $J$.  Since
the conclusions hold for $\varphi_1$ and $\varphi^\dagger_1$ by
\cite[p.~76]{Jac:J3} and $\times$ is bilinear, we are done.
\end{pf}

\section{Brown algebras and groups of type $E_6$}

\begin{defn} \cite[p.~135]{A:structintro}, \cite[1.1]{AF:CD}
Suppose that $(A, {-})$ is a finite-dimensional (and perhaps
nonassociative) $F$-algebra with $F$-linear involution.  For $x, y \in
A$, define  
$V_{x,y} \in \EndF(A)$ by 
\begin{neqn} \label{braceprod}
V_{x,y} z := 
\{ x, y, z \} := (x\overline{y})z + (z\overline{y})x -
(z\overline{x})y,
\end{neqn}
for $z \in A$.  One says that $(A, {-})$ is a {\em structurable algebra}
if 
\[
[V_{x,y}, V_{z,w}] = V_{V_{x,y}z, w} - V_{z,V_{y,x}w}.
\]

The multiplication algebra of $(A, {-})$ is the (associative)
subalgebra of $\EndF(A)$ generated by the involution $-$, left
multiplications by elements of $A$, and right multiplications by
elements of $A$.  If the center of the multiplication algebra of $(A,
{-})$ is $F$, then $(A, {-})$ is said to be {\em central}.

We say that $(A, {-})$ is {\em simple} if it has no two-sided ideals
which are stabilized by $-$.
\end{defn}

This definition of a structurable algebra in terms of this $V$
operator may seem unmotivated.  There is, however, an
alternative (partial) characterization which works for the case that
we are interested in.  Suppose that $(A, {-})$ is an $F$-algebra with
$F$-linear involution which is generated as an $F$-algebra by its
space of symmetric elements.  Then by \cite[p.~144]{A:structintro}
$(A, {-})$ is structurable if and 
only if it is skew-alternative (i.e., $[s, x, y] = -[x, s, y]$ for all
$x, y \in A$ and $s$ skew in $A$ where $[x,y,z] := (xy)z - x(yz)$) and
it supports a symmetric bilinear form $\qform{,}$ which satisfies
\[
\qform{\bar{x}, \bar{y}} = \qform{x,y} \hbox{\ and\ } \qform{zx,y} =
\qform{x, \bar{z}y}
\]
for all $x$, $y$, $z \in A$.

Basic examples of structurable algebras are Jordan algebras (with
involution the identity) and central simple algebras with involution.
For Jordan algebras, the ternary product $\{ , , \}$ given in
(\ref{braceprod}) is the usual triple product as in
(\ref{jordisodef}) or \cite[p.~36, (58)]{Jac:J} and the symmetric
bilinear form is  the trace form $T$.

\begin{eg} \cite[1.9]{A:skew} \label{basic}
Let $J$ be an Albert $F$-algebra and $\zeta \in \Fx$.
We define a structurable algebra $(B, {-}) := \B(J, F \times F,
\zeta)$ by setting 
$B$ to be the vector space 
\[
\tbtmat{F}{J}{J}{F}
\]
with multiplication given by
\[
\tbtmat{\alpha_1}{j_1}{j'_1}{\beta_1}
\tbtmat{\alpha_2}{j_2}{j'_2}{\beta_2} = 
\tbtmat{\alpha_1 \alpha_2 + \zeta T(j_1, j'_2)}{\alpha_1 j_2 + \beta_2 j_1 +
\zeta (j'_1 \times j'_2)}{\alpha_2 j'_1 + \beta_1 j'_2 + j_1 \times
j_2}{\beta_1 \beta_2 + \zeta T(j_2, j'_1)}.
\]
Endow $B$ with the involution
$-$ given by
\[
\overline{\tbtmat{\alpha}{j}{j'}{\beta}} =
\tbtmat{\beta}{j}{j'}{\alpha}.
\]
We use the abbreviation $\B(J, F \times F)$ for $\B(J, F \times F, 1)$.
\end{eg}

This is a central simple structurable algebra, and is denoted by
$\mathcal{M}(\zeta T, \zeta N, \zeta^2 N)$ in \cite{AF:CD} for $T$ and
$N$ the trace and norm of $J$ respectively.

The study of such
algebras precedes the notion of structurable algebras significantly:
structurable algebras were introduced in \cite{A:structintro} and these
algebras are a special case of those discussed in \cite{Brown:th} and
\cite{Brown:newalg}.  To be precise, the algebras that Brown studied
involved parameters $\mu$, $\nu$, $\omega_1$, $\omega_2$, $\delta_1$,
and $\delta_2$.  If one sets $\mu = \nu = \omega_1 = 1$ and $\omega_2 = \delta_1 = \delta_2 = \zeta$, the algebra
$\B(J, F \times F, \zeta)$ is obtained.

\begin{eg} \label{basic2}
If $\D$ is a quadratic field extension of $F$, we define a
structurable algebra $\B(J, \D)$. 
There is an ``outer'' automorphism $\varpi$ of $\B(J, F \times F)$ given by
\begin{neqn} \label{outeraut} \label{flipaut}
\varpi \basmat := \tbtmat{\beta}{j'}{j}{\alpha}.
\end{neqn}
Let $\iota$ denote the unique nontrivial
$F$-automorphism of $\D$ and set $\B(J, \D)$ to be the $F$-subalgebra of the
$\D$-algebra $\B(J, F \times F) \otimes_F \D$ fixed by $\varpi \otimes \iota$.
Then $\B(J, \D)$ is a structurable algebra over $F$ and
\begin{neqn}
\B(J, \D) \otimes_F \D \cong \B(J \otimes_F \D, \D \times \D).
\end{neqn}
\end{eg}

Thus the algebra $\B(J, \D)$ is also a central
simple structurable algebra by descent.

\begin{defn} \label{browndef}
For $J^d$ the split Albert algebra over $F$, we call $\B^d := \B(J^d,
F \times F)$ the {\em split Brown algebra} over $F$ and $\B^q_{\D/F}
:= \B(J^d, \D)$ the {\em quasi-split Brown algebra with inner
extension $\D$}. 

We say that an $F$-algebra with involution $(B, {-})$ is a {\em Brown
algebra} if $(B, {-}) \otimes_F F_s \cong 
\B^d \otimes_F  F_s$ for $F_s$
a separable closure of $F$.
\end{defn}

By the classification of central simple structurable algebras due to
Smirnov (see \cite{Smirnov:AeL} or \cite{Smirnov:CM}) and Allison, if $F$ has
characteristic $\ne 5$ (and, as always, $\ne 2, 3$), we could have equally well
defined a Brown $F$-algebra to be a central simple structurable
algebra over $F$ of dimension 56 and skew-dimension 1 (i.e.~the space
of skew-symmetric elements is 1-dimensional). I do not know of a
classification theorem for central simple structurable algebras in characteristic 5.

In any event, any Brown algebra $\bbrown$ has a 1-dimensional space of
skew-symmetric elements.  If $s_0 \in B$ spans this space, then $s_0^2 \in
\Fx$, which one can see by descent from the split case or see
\cite[2.1(b)]{AF:CD} for a different argument.  We say that $\bbrown$
is {\em of type   
  1} if $s_0^2$ is a square in $F$ and that it is of {\em of type 2}
otherwise.  We call $\D := F[s_0]$ the {\em discriminant algebra} of
$\bbrown$. 

It is worth mentioning that not all Brown algebras of type
2 are as in Example \ref{basic2}, see \ref{nonredrmk}. 

\begin{lem} \label{isoprop} \label{e6lem} 
\begin{enumerate}
\item {\rm \cite[4.5]{AF:CD}} Any Brown algebra of type 1 is isomorphic to some algebra of the
form $\B(J, F \times F, \zeta)$.

\item For Albert $F$-algebras $J_1$ and $J_2$, $\B(J_1, F\times F,
  \zeta_1) \cong \B(J_2, F \times F, \zeta_2)$ if and only if there is
  a norm similarity $J_1 \longto J_2$ which has multiplier $\zeta_1 /
  \zeta_2$ or $\zeta_1 / \zeta^2_2$.

\item If $J$ is reduced then $\B(J, F \times F, \zeta)
\cong \B(J, F \times F)$ for all $\zeta \in \Fx$.
\end{enumerate}
\end{lem}

Before we proceed with the proof, it is should be noted that any obvious
analogue of statement (2) for Brown algebras of type 2
is false.  Specifically, let $J_1 :=
\h_3(\C, 1)$ and $J_2 := \h_3(\C, \gamma)$ for $\C$ the Cayley
division algebra over the reals and $\gamma := \diag (1, -1,
1)$.  Then $J_1$ and $J_2$ have isometric norms.  However, for $\B_i
:= \B(J_i, \mathbb{C})$, the Lie algebra of $\autp (\B_i)$ is classically
denoted by $\L(J_i)_{-1}$.  The two Lie algebras $\L(J_i)_{-1}$ for $i
= 1,2$ are not isomorphic since they have different signatures, see
\cite[pp.~119, 120]{Jac:ex}.
\medskip

\begin{pf}
(1) boils down to the main results of the paper \cite{Sp:cubic} of
Springer, see \cite[4.5]{AF:CD}.

\smallskip
(3): We have an isomorphism $\pi \!: \B(J, F \times F, \zeta) \longto
\B(J, F \times F, \zeta^2)$ given by

\smallskip
(2), $\Implies$: Suppose that $\tau \!: \B(J_1, F \times F, \zeta_1) \longto
\B(J_2, F \times F, \zeta_2)$ is an algebra isomorphism.  We may think of 
\[
s_0 := \tbtmat{1}{0}{0}{-1}
\]
as a skew-symmetric element in $\B(J_i, F \times F)$ for $i = 1, 2$.
Since $\tau(s_0)$ must also be skew-symmetric and $\tau(s_0)^2 =
\tau(s_0^2) = 1$, we must have that $\tau(s_0) = \pm s_0$.  

Suppose first
that $\tau(s_0) = s_0$.  Then $\tau$ fixes the diagonal matrices
elementwise and 
so $\tau$ is given by
\[
\tau \basmat = \tbtmat{\alpha}{\varphi(j)}{\varphi'(j')}{\beta}
\]
for some linear maps $\varphi$, $\varphi' \!: J_1 \longto J_2$.  Since
$\tau$ is an algebra isomorphism, $\zeta_2 T_2(\varphi(j_1), \varphi'(j'_1)) =
\zeta_1 T_1(j_1, j'_1)$ for $T_i$ the trace form on $J_i$ and $j_1$, $j'_1 \in
J_1$.  Let $\varphi^\dagger \!: J_1 \longto J_2$ be the unique linear
map such that $T_2(\varphi(j_1), \varphi^\dagger(j'_1)) = T_1(j_1, j'_1)$
for all $j_1, j'_1 \in J_1$.  Then $\varphi' = (\zeta_2 / \zeta_1)
\varphi^\dagger$.  Since
\[
\varphi(j_1) \times \varphi(j'_1) = \frac{\zeta_1}{\zeta_2}
\varphi^\dagger (j_1 \times j'_1),
\]
we have 
\[
N_2(\varphi(j)) = \frac{\zeta_1}{\zeta_2} N_1(j).
\]

Now suppose that $\tau(s_0) = -s_0$.  We have an isomorphism 
$\pi \!: \B(J_2, F \times F, \zeta_2) \iso \B(J_2, F \times F,
\zeta^2_2)$ given by 
\[
\pi\basmat = \tbtmat{\beta}{j'}{\zeta_2^{-1} j}{\alpha}.
\]
Composing $\tau$ with $\pi$, we get a new isomorphism
$\pi \tau \!: \B(J_1, F \times F, \zeta_1) \longto \B(J_2, F \times F,
\zeta^2_2)$  
such that $\pi \tau(s_0) = s_0$.  By our previous result, there must
be a norm similarity between $J_1$ and $J_2$ with multiplier $\zeta_1 /
\zeta^2_2$. 

\smallskip
(2), $\If$: Conversely, if one is given such a norm similarity
$\varphi$, then one can run the argument just given for the other
direction backwards to produce the desired algebra isomorphism.

\smallskip
(3) is clear from (2), since for any $\zeta \in \Fx$, $J$ has a norm
similarity with multiplier $\zeta$ as in (\ref{jordmult}).
\end{pf}

\begin{thm} \label{e6thm}
\begin{enumerate}
\item If $\bbrown$ is a Brown algebra of type $t$ over $F$, then
  $\autp\bbrown$ is a 
simply connected group of type $^tE_6$ over $F$ with trivial Tits
algebras. Every 
simply connected group of type $E_6$ with trivial Tits algebras arises
in this way. 

\item The automorphism group of the
split Brown algebra, $\autp(\B^d)$, is the split simply connected
group of type $E_6$.   For $\D$ a quadratic field extension of $F$,
$\autp(\B^q_{\D/F})$ is the quasi-split simply connected group of
type $^2E_6$ with inner extension $\D$.
\end{enumerate}
\end{thm}

\begin{pf}
(2): Let $f$ be an $F$-algebra automorphism of $\B^d$.  Since $f$
respects the involution ${-}$ on $\B^d$, it must map the 1-dimensional
subspace $\sym$ of skew-symmetric elements to itself.  Since $\sym$
and $F \cdot 1$ span the diagonal elements of $\B^d$, $f$ must
preserve these.  Set $u := \sozmat$.  Then the diagonal matrices are
spanned by $u$ and 
$\ubar$.  Set $\stbtmat{a}{0}{0}{b} := f(u)$.  Then
\[
\tbtmat{ab}{0}{0}{ab} = f(u) \overline{f(u)} = f(u) f(\ubar) = 0.
\]
Thus $a = 0$ or $b = 0$, but not both.

Suppose for the moment that $b = 0$.  Then since $f$ is a Brown
algebra automorphism, 
\begin{neqn} \label{firstf}
f\stbtmat{\alpha}{j}{j'}{\beta} =
\stbtmat{\alpha}{\varphi(j)}{\varphi^\dagger(j')}{\beta}
\end{neqn}
for $\varphi$
some norm isometry of $J^d$.

Otherwise $a = 0$ and $\varpi f$ is of the form just described in the
``$b = 0$'' case.

Thus $\autp(B^d)$ is isomorphic to a semidirect product of $\Z_2$ and
the group of norm isometries of $J^d$, which is known to be split
simply connected of type $E_6$. 

Now consider $\B^q := \B(J^d, \Delta)$.  Since $\B^q \otimes_F \Delta
\cong \B^d \otimes_F \Delta$, $\autp(\B^q)$ is a simply connected
of type $E_6$.  It has trivial Tits algebras since any Tits algebra
would have 
exponent a power of 3 and would be split by $\Delta$.  The center of
$\autp(\B^d \otimes_F \Delta)$ consists of the maps $f_\omega$, where
\[
f_\omega \basmat = \tbtmat{\alpha}{\omega
j}{\omega^2 j'}{\beta}
\]
for $\omega$ a cube root of unity.  To see that $\autp(\B^q)$ is of
type 2, we need to see that the Galois action on the center of
$\autp(B^d \otimes_F \Delta)$ is not the same as the action induced by the
$\iota$-semilinear automorphism $\varpi \otimes \iota$ of $\B^d
\otimes_F \Delta$ 
which defines $\B^q$:
\[
(\Id \otimes \iota) f_\omega (\Id \otimes \iota) \basmat \ne
\tbtmat{\alpha}{\iota(\omega)^2 j}{\iota(\omega) j'}{\beta} =
(\varpi \otimes \iota) f_\omega (\varpi \otimes \iota) \basmat \ne
\]
This shows that $\autp(\B^q)$ is simply connected and Tits-trivial of
type $^2E_6$ with inner extension $\Delta$.

We have an injection $\aut(J^d) \injects \autp(\B^q)$ given by
$\varphi \mapsto f$ as in (\ref{firstf}) which produces a
rank four $F$-split torus in $\autp(\B^q)$.  The Galois group of
$\Delta$ over $F$ acts nontrivially on the set of simple roots of the
Dynkin diagram of $\autp(\B^q)$ which is of type $E_6$, so it has precisely
four orbits.  Thus $\autp(\B^q)$ is quasi-split 
\cite[Lem.~4.2]{G:iso}. 

\smallskip
(1):  Since $\autp\bbrown$ is a form of $\autp(\B^d)$ which is simply
connected of type $E_6$, so is $\autp\bbrown$.  Suppose that $\Delta$
is the discriminant algebra of $\bbrown$ so that $\D = F \times F$ if
$\bbrown$ is of type 1 and $\bbrown \otimes_F \D$ is of type 1
otherwise.  Let $K$ be a field extension of $F$ which generically
quasi-splits $\autp\bbrown$ as in \cite{KR}.  Then $\autp\bbrown
\times_F K \cong \autp (\B(J^d \otimes_F K, \Delta\otimes_F K))$,
which verifies the type of $\autp\bbrown$ since $F$ is algebraically
closed in $K$.  That  one obtains every group of type $E_6$ in this
manner follows by the usual Galois cohomology argument.
\end{pf}

\begin{rmk}
Zinovy Reichstein 
suggested to me that there should be an 
invariant
\[
g_4 \!: H^1(F, G) \longto H^4(F, \mmu3)
\]
defined for $G$ simply connected of type $E_6$.  The theorem we just
proved allows us to sketch the definition of such an invariant here in
the case where $G$ is quasi-split. 
By \cite[12.13]{Rei:ed}, this shows that the 3-primary essential
dimensional (see \cite[3.1]{Rei:ed} for a definition) of a simply
connected group of type $E_6$ is $\ge 4$, hence that the overall
essential dimension of such a group is $\ge 4$.
(This result has also been obtained by
Reichstein and Youssin by other means, see \cite[8.19.4]{RY}.)

We first observe that as in the construction of the Serre-Rost
invariant
\begin{neqn} \label{SRinv}
g_3 \!: H^1(F, \aut(J^d)) \longto H^3(F, \Z / 3)
\end{neqn}
in \cite{Ro} or \cite{PR:el}, since 2 and 3 are coprime we can extend
scalars up to a quadratic extension, define $g_4$ up there, and then
use the corestriction to define $g_4$ over the ground field.  In
particular, this reduces 
us to considering the case where $G$ is actually split.  By Theorem
\ref{e6thm}(2), $H^1(F, G)$ classifies pairs
$(\B, \phi)$ where $\B$ is a Brown $F$-algebra of type 1 and $\phi$ is 
an algebra isomorphism from the discriminant algebra of $\B$
to $F \times F$.  By Lemma \ref{e6lem}, we can write any such
algebra $\B$ as $\B(J, F \times F, \zeta)$ where $\phi$ is given by
$\stbtmat{\alpha}{0}{0}{\beta} \mapsto (\alpha, \beta)$.  Then we
define
\[
g_4(\B, \phi) := g_3(J) \cup (\zeta),
\]
where $(\zeta) \in H^1(F, \mmu3)$ and we are identifying $\Z / 3
\otimes \mmu3 = \mmu3$.  

We must check that this is well-defined.
Suppose that $(\B, \phi) \cong (\B, \phi')$ where $\B' \cong \B(J', F
\times F, \zeta')$.  Then by Lemma \ref{e6lem}(2), there is a norm
similarity $J \longto J'$ with multiplier $\zeta / \zeta'$.  (The
other possibility of having multiplier $\zeta / (\zeta')^2$ doesn't
occur since this requires a switch in the identification of the
diagonal matrices with $F \times F$.)  Again we may extend scalars to
a quadratic extension so that $J$ is a first Tits construction (a
quadratic extension is sufficient by, for example,
\cite[39.19]{KMRT}).  Then since $J$ and $J'$ are isotopic, they
are even isomorphic by \cite[4.9]{PR:sp}, so
\begin{neqn} \label{invdiff}
\left( g_3(J) \cup (\zeta) \right) - \left( g_3(J') \cup (\zeta')
\right) = g_3(J) \cup (\zeta / \zeta').
\end{neqn}
However, $J$ has a norm similarity with multiplier $\zeta / \zeta'$,
so its norm represents $\zeta / \zeta'$.  By \cite[4.6]{PR:sp}, $J
\cong J(A, \lambda)$ for some central simple $F$-algebra $A$ of degree 
3 and some $\lambda \in \Fx$ such that $\zeta / \zeta' \in
\NrdA(\Ax)$, so $(A) \cup (\zeta / \zeta') = 0$.  
Since $g_3(J) = (A) \cup (\lambda)$, it follows that the difference in
(\ref{invdiff}) of our two possibilities is zero, so our map $g_4$ is
well-defined. 
\end{rmk}

\section{Background on Freudenthal triple systems}

Multiple authors have studied Freudenthal triple systems as a means to
understanding groups of type $E_7$.  Axiomatic treatments appear in
\cite{Brown:E7}, \cite{Meyb:FT}, and \cite{Ferr:strict}, for example.
These authors considered a general sort of Freudenthal triple system,
but we are only interested in a particular kind.

\begin{defn} (Cf.~\cite[p.~314]{Ferr:strict})
A {\em (simple) Freudenthal triple system} is a 3-tuple $(V, b,
t)$ such that $V$ is a 56-dimensional vector space, $b$ is a
nondegenerate skew-symmetric bilinear form on $V$, and $t$ is a
trilinear product $t \!: V \times V \times V \longto V$.

We define a 4-linear form $q(x,y,z,w) := b(x, t(y, z, w))$ for $x$, $y$, $z$, $w \in V$, and we 
require that
\begin{description}
\item[FTS1] $q$ is symmetric,
\item[FTS2] $q$ is not identically zero, and
\item[FTS3] $t(t(x,x,x),x,y) = b(y,x) t(x,x,x) + q(y, x,x,x) x$ for
all $x, y \in V$.
\end{description}
We say that such a triple system is {\em nondegenerate} if the quartic
form $v \mapsto q(v,v,v,v)$ on $V$ is absolutely irreducible (i.e.,
irreducible over a separable closure of the base field) and {\em
degenerate} otherwise.
\end{defn}

Note that since $b$ is nondegenerate, FTS1 implies that $t$ is symmetric.

\begin{eg}
For $J$ an Albert $F$-algebra and $\zeta \in \Fx$, we can construct a
Freudenthal triple system as follows.  Set 
\[
V := \tbtmat{F}{J}{J}{F}.
\]
As in \cite[1.10]{A:skew} or
\cite[(5), (6), p.~87]{Brown:E7}, for  
\[
x_1 = \tbtmat{\alpha_1}{j_1}{j'_1}{\beta_1} \hbox{\ and\ }
x_2 = \tbtmat{\alpha_2}{j_2}{j'_2}{\beta_2},
\]
set
\begin{neqn} \label{traceform}
b(x_1, x_2) = (\alpha_1 \beta_2 - \alpha_2 \beta_1) + \zeta (T(j_1,
j'_2) - T(j'_1, j_2)).
\end{neqn}
Let $q$ be as in the definition of a Freudenthal triple system.  Since
$b$ is nondegenerate, if we 
know $b$ and $q$ then we can determine $t$, at least in principle.  In
this case, for
\[
x = \tbtmat{\alpha}{j}{j'}{\beta},
\]
we define
\begin{neqn} \label{normform}
q(x,x,x,x) = 12\left(4 \alpha \zeta N(j) + 4 \beta \zeta^2 N(j') - 4\zeta^2
T({j'}^{\#}, j^{\#}) + (\alpha \beta - \zeta T(j, j'))^2\right).
\end{neqn}

Then $(V,b,t)$ is a nondegenerate Freudenthal triple system.
When $\zeta = 1$, it is denoted by $\M(J)$.
\end{eg}

By \cite[\S 4]{Brown:E7}, a Freudenthal triple system is nondegenerate
if and only if it is a form of a triple system $\M(J)$ for some $J$
(i.e., it becomes isomorphic to $\M(J)$ when one extends scalars to a
separable closure).  There is a way to distinguish between the two
sorts of triple systems over the ground field, see \cite[\S 2]{G:e7}.
Meyberg \cite[\S 7]{Meyb:FT} uses a different terminology; he says
that the the 
nondegenerate ones are ``of main type''.
It is the nondegenerate ones which are relevant to groups of type
$E_7$.  (It is not clear what the automorphism group of a degenerate
triple system is, but
it is not connected and it contains --- at least over a separably
closed field --- a torus of rank 28, see \cite[\S 1]{G:e7}.)

For $\M$ a Freudenthal triple system over $F$, we follow
Freudenthal's notation \cite[\S 3]{Frd:E7.1} and write $\Inv(\M)$ for
the algebraic group whose $F$-points are the isomorphisms of $\M$.

\begin{thm} \label{FTSthm}
For $\M$ a nondegenerate Freudenthal triple system over $F$,
$\Inv(\M)$ is a simple simply connected algebraic
group of type $E_7$ with trivial Tits algebras.  This construction
produces all such groups.  Moreover, $\Inv(\M(J^d))$ is split.
\end{thm}

\begin{pf}
The group $\Inv(\M)$ is simple by \cite[p.~100,
Thm.~6]{Brown:E7} and it is of type $E_7$ by \cite[\S 5]{Frd:E7.1}.
Since it has center $\mmu2$, it is simply connected.  A simply
connected group of type $E_7$ is isomorphic to $\Inv(\M)$ for some
$\M$ if and only if it has trivial Tits algebras by
\cite[6.5.2]{Ti:R}.

Finally we show that $\Inv(\M^d)$ is split for $\M^d := \M(J^d)$.  For
an arbitrary Albert $F$-algebra $J$ and a norm similarity $\varphi$ of
$J$ with multiplier $\lambda \in \Fx$, we have an element $f_\varphi
\in \Inv(\M(J))$ given by
\begin{neqn} \label{Ts} \label{fdef}
f_\varphi \basmat := \tbtmat{\lambda^{-1}
  \alpha}{\varphi(j)}{\varphi^{\dagger}(j')}{\lambda \beta}.
\end{neqn}
This restricts to an injection $\Inv(J) \injects \Inv(\M(J))$.

Now, $\Inv(J^d)$ is split of type $E_6$, so let $S_6$ denote the image
in $\Inv(\M^d)$ of a rank 6 split torus in $\Inv(J^d)$.  Let $S_1$ be
the image of $\GmF$ in $\Inv(\M^d)$ under the map $x \mapsto f_{L_x}$
where $L_x$ denotes left multiplication by $x$.  Then $S := S_1 S_6$
is a rank 7 $F$-split torus in $\Inv(\M^d)$, so the group is split.
\end{pf}

Another popular approach to constructing groups of type $E_7$ is by
realizing the groups as automorphism groups of a particular quartic
form which is not the one used here, see \cite[\S\S 7, 8]{Asch:mult}
and \cite{Coop:E7}.  The precise relationship between their approach
and ours is not clear.

\section{Brown algebras and groups of type $E_7$} \label{e7sec}

Let $\bbrown$ be a Brown algebra over $F$.  By definition, the space
$S$ of skew-symmetric elements of $B$ has dimension 1, so pick some
nonzero $s_0 \in B$ such that $S = Fs_0$.  There is a natural map
$\psi \!: B \times B \longto S$ given by
\begin{neqn} \label{allpsidef}
\psi(x,y) := x\overline{y} - y\overline{x}.
\end{neqn}
It is known that for any $x \in B$, 
\begin{neqn} \label{skewass}
s_0 (s_0 x) = (x s_0) s_0 = x(s_0^2) \hbox{\ and\ } s_0^2 = \mu
\end{neqn}
for some $\mu \in \Fx$ by descent or by \cite[p.~135,
Prop.~1]{A:structintro}.  Thus
the map
\begin{neqn} \label{bform}
b(x,y) := \psi(x,y) s_0
\end{neqn}
provides a skew-symmetric bilinear form on $B$.
We also have a trilinear map $t \!: B \times B
\times B \longto B$ given by 
\begin{neqn} \label{tform}
t(y, z, w) := 2\{ y, s_0 z, w\} - b(z,w) y - b(z,y) w - b(y,w)z.
\end{neqn}
Then $b$ and $t$ give $B$ the structure of a (Freudenthal) triple
system by \cite[2.18]{AF:CD}, and we say that $(B,b,t)$ is a triple
system {\em associated to} $\bbrown$.  Of course, this triple system is not
uniquely determined: if 
we choose some other element $s'_0 \in S$ such that $S = Fs'_0$, then
$s'_0 = \lambda s_0$ for some $\lambda \in \Fx$, and this choice of
$s'_0$ would give us a skew-symmetric bilinear form $b'$ and a
trilinear form $t'$ such that $b = \lambda b'$ and $t = \lambda t'$.

The definitions of $b$ and $t$ given above may appear to be ad hoc,
but in fact they are not.  For any structurable algebra with
involution $-$, there is a
natural symmetric bilinear trace form given by setting $\qform{x, y}$
to be the trace of left multiplication by $x\overline{y} + y
\overline{x}$.  This trace form was instrumental in the classification
of central simple structurable algebras in \cite{A:structintro}.  In the
special case where the structurable algebra is a Brown algebra, any
nonzero skew-symmetric element $s_0$ spans the skew elements, and $b$
is a scalar multiple of the map $(x, y) \mapsto \qform{s_0 x, y}$.
There is also a norm form on any central simple structurable algebra
defined up to a scalar \cite{AF:norm}.  This norm specializes to the
usual norm in the case where the structurable algebra is Jordan and to
the reduced norm when it is an associative central simple algebra.  If
we write $\nu$ for this norm on a Brown algebra, then by
\cite[2.17]{AF:CD} the trilinear map $t$
is given by
\[
b(x, t(x,x,x)) = 12 \mu \nu(x)
\]
where $\mu := s^2_0 \in \Fx$ as above.

\begin{eg}
Consider the Brown algebra $\B := \B(J, F \times F, \zeta)$ and let
$s_0 = \stbtmat{1}{0}{0}{-1}$.  Then a triple system associated to $\B$
is given by the formulas (\ref{traceform}) and (\ref{normform})
\cite[pp.~192--195]{AF:CD}.  Since this triple system is a form of
$\M(J)$ and all Brown algebras are a form of $\B$, every Freudenthal
triple system associated to a Brown algebra is nondegenerate.
\end{eg}

It is these triple systems which
will provide us with our groups of type $E_7$, and so we need to
understand when two Brown algebras provide us with the same group.

\begin{defn}\label{simdef}
Let $(B_1, {-})$ and $(B_2, {-})$ be two Brown algebras over $F$, and
let $(B_i, b_i, t_i)$ be a triple system related to $B_i$ for $i = 1,2$.
We say that $(B_1, {-})$ and $(B_2, {-})$ are {\em similar} if there
is some $F$-vector space map $f \!: B_1 \longto B_2$ and some $\lambda
\in \Fx$ such that 
\begin{neqn} \label{simeq}
b_2(f(x), f(y)) = \lambda b_1(x,y) \hbox{\ and\ } t_2(f(x), f(y),
f(z)) = \lambda f(t_1(x,y,z))
\end{neqn}
for all $x$, $y$, $z \in B_1$.  Such an $f$ is said to be a {\em
similarity}.  

If $(B_1, {-}) = (B_2, {-})$, $(B_1, b_1, t_1) = (B_2, b_2, t_2)$,
and $f$ satisfies 
(\ref{simeq}), then we say that $f$ is a {\em similarity of $(B_i,
  {-})$ with multiplier $\lambda$}.  If $\lambda = 1$, $f$ is called
an {\em isometry}.  
\end{defn}

Note that these definitions are independent of the choice of triple
systems for our Brown algebras.

We define three algebraic groups associated to a Brown algebra $\B$.
Set $\GInv(\B)$ to be the algebraic group with $F$-points
\[
\GInv(\B)(F) = \{ f \in \EndF(\B) \mid \hbox{$f$ is a similarity} \}
\]
and let $\Inv(\B)$ be the algebraic group with $F$-points
\[
\Inv(\B)(F) = \{ f \in \EndF(\B) \mid \hbox{$f$ is an isometry} \}.
\]
Set $\PInv(\B)$ to be the quotient of $\GInv(\B)$ by its center.


Two Albert algebras have similar norms if and only if they are
isotopic.  There is a natural notion of isotopy for structurable
algebras, and we recall the definition here.

\begin{defn} \cite[p.~132]{AH}
Let $(A, {-})$ and $(A', {-})$ be two structurable algebras.  They are said to be 
{\em isotopic} (abbreviated $(A, {-}) \sim (A', {-})$)
if there are $F$-linear maps $\alpha, \beta \!: A \longto A'$ such that
\begin{neqn} \label{isodefeq}
\alpha \{ x, y, z \} = \{ \alpha(x), \beta(y), \alpha(z) \}'.
\end{neqn}
\end{defn}

In the special case where $A = A'$ as vector spaces and the two
structurable algebras are actually Jordan, this reduces to the
standard notion of isotopy.  One direction is clear enough:
substituting $y = 1$ in (\ref{isodefeq}), we see that this definition
specializes to that given in (\ref{jordisodef}) with $u = \beta(1)$.
Conversely, any isotopy as classically defined for Jordan
algebras induces an isotopy as we have just defined them, see
\cite[p.~83]{AH}.

It turns out that the notions of similarity and isotopy are equivalent
for Brown algebras, so that we have a 
situation analogous to that with Albert algebras.

\begin{prop} \label{isosim}
Two Brown algebras are similar if and only if they are isotopic.
\end{prop}

\begin{pf}
Let $\B := (B, {-})$ and $\B' := (B', {-}')$ denote our two Brown algebras, and
let $(B,b,t)$ and $(B',b',t')$ provide triple systems associated to them.

( $\Implies$ ): Let $f \!: B \longto B'$ and $\lambda$ be as in
Definition \ref{simdef}.  Let $s_0$, $s'_0$ span the skew-symmetric
elements in $\B$ and $\B'$.  We may assume that $b$, $t$ and $b'$,
$t'$ are given by the formulas in (\ref{bform}) and (\ref{tform}).
Then plugging this into (\ref{simeq}) provides 
\[
\{ f(y), s'_0 f(z), f(w) \}' = \lambda f \{ y, s_0 z, w\}
\]
for all $x$, $y$, $z \in B$.

We replace $z$ by $(\lambda s_0)^{-1} z$ and $f$ by $g$ defined by
$g(z) := s'_0 f( (\lambda s_0)^{-1} z)$ 
to obtain
\[
\{ f(y), g(z), f(w) \}' = f \{ y, z, w \}.
\]
Thus $f$ is an isotopy.

\newcommand{\bconj}{\B^{\qform{u}}}
( $\If$ ): Since $\B$ and $\B'$ are isotopic, there is
some element $u \in \B$ such that $\B'$ is isomorphic to a new Brown
algebra denoted by $\bconj$ \cite[p.~134, Prop.~8.5]{AH}.  This algebra is a 
Brown algebra over $F$ with the same underlying vector space as
$\B$.  Fix some $s_0$ which spans the space of skew-symmetric
elements of $B$.  Then $s_0^{\qform{u}} := s_0 u \ne 0$ spans the
space of skew-symmetric elements of $\bconj$ by \cite[1.16]{AF:CD}.
Let $(B,b,t)$ and $(B, b^{\qform{u}}, t^{\qform{u}})$ be the triple
systems associated to $\B$ and $\bconj$ determined by $s_0$ and
$s^{\qform{u}}_0$ by the formulas (\ref{bform}) and (\ref{tform}).

Let $\psi^{\qform{u}}$ be a map on $\bconj$ defined as
in (\ref{allpsidef}), and let $\nu$ denote the conjugate norm on $\B$, so
$\nu(x) = \frac{1}{12\mu} q(x, x, x, x)$ for $\mu := s^2_0 \in \Fx$ by
\cite[2.17]{AF:CD}.  (Note that this definition of $\nu$ is
independent of the choice of $s_0$.) For $x, y \in B$, let $\lambda
\in \Fx$ be such that $\psi(x,y) = \lambda s_0$. 
Then
\begin{neqn} \label{isoteq}
\begin{array}{rclr}
b^{\qform{u}} (x,y) &=& \psi^{\qform{u}}(x,y) s^{\qform{u}}_0 &\\
&=& (\psi(x,y)u)s^{\qform{u}}_0&\hfill\hbox{by \cite[1.17]{AF:CD}}\\
&=& \lambda (s_0^{\qform{u}})^2 \\
&=& \lambda \nu(u) \mu & \hbox{by \cite[3.2]{AF:CD}}\\
&=& \nu(u) b(x,y).
\end{array}
\end{neqn}
Let $\mu^{\qform{u}} := (s^{\qform{u}}_0)^2 = \nu(u) \mu$
and let $\nu^{\qform{u}}$ be the conjugate norm on $\bconj$.  By
\cite[3.7]{AF:CD}, $\nu^{\qform{u}}(x) = \nu(u) \nu(x)$ for all $x \in
B$.  We can
linearize $\nu^{\qform{u}}$ to get a unique symmetric 4-linear form
such that $\nu^{\qform{u}}(x,x,x,x) = 24 \nu^{\qform{u}}(x)$.  Then
\[
\begin{array}{rclr}
b^{\qform{u}} (x, t^{\qform{u}} (y, z, w)) &=&
\frac{\mu^{\qform{u}}}2 \nu^{\qform{u}} (x,y,z,w) &
\\
&=& \frac{\nu(u)^2 \mu}2 \nu(x,y,z,w)&
\\
&=& \nu(u)^2 b(x, t(y,z,w)).&
\end{array}
\]
By the nondegeneracy of $b$ and (\ref{isoteq}), the identity map on
$\B$ is a similarity with multiplier $\nu(u)$.
\end{pf}

\begin{lem} \label{e7lem}
\begin{enumerate}
\item $\B(J, F \times F, \zeta) \sim \B(J, F \times F)$.

\item $\B(J_1, F \times F) \sim \B(J_2, F \times F)$ if and only if $J_1 \sim
J_2$.
\end{enumerate}
\end{lem}

\begin{pf}
(1): Define an $F$-vector space isomorphism $f \!: \B(J, F \times F) \longto \B(J, F \times F, \zeta)$ via
\begin{neqn} \label{fsim}
f\tbtmat{\alpha}{j}{j'}{\beta} = \tbtmat{\zeta \alpha}{j}{j'}{\beta}.
\end{neqn}
A quick check of formulas (\ref{traceform}) and (\ref{normform}) show
that this is a similarity of the triple systems.

(2): One of the triple systems associated to $\B(J_i, F \times F)$ is
$\M(J_i)$.  Thus $\B(J, F \times F) \sim \B(J', F \times F)$ if and
only if $\M(J_1)$ is similar to $\M(J_2)$.  By
\cite[6.8]{Ferr:strict}, this occurs if and only if $J_1 \sim J_2$. 
\end{pf}

The construction of Freudenthal triple systems from Brown algebras is
already interesting, and the following lemma makes it a useful tool.

\begin{lem} \label{goodconst}
Every Freudenthal triple system is associated to some Brown algebra.
\end{lem}

\begin{pf}
Let $\M = (V,b,t)$ be a Freudenthal triple system.  If $\M \cong
\M(J)$ for some Albert $F$-algebra  
$J$, then $\M$ is obtained from $\B(J, F \times F)$ and we are done.
Otherwise, we pick some $v \in V$ such that $q(v) \ne 0$.  Since $\M
\not\cong \M(J)$ for any $J$, $d := q(v) / 12$ is not a square in
$\Fx$ by \cite[3.4]{Ferr:strict}.  Set $\D := F(\sqrt{d})$.  Then $\M
\otimes_F \D \cong \M(J)$ for some Albert $\D$-algebra $J$.  Moreover,
for $\iota$ the nontrivial  
$F$-automorphism of $\D$, $\M$ is $F$-isomorphic to the fixed points of $\M(J)$
under the $\iota$-semilinear automorphism given by 
\[
\eta \basmat = 
\tbtmat{\iota(\beta/\delta)}{\varphi^\dagger
  \iota(j')}{-\varphi\iota(j)}{\delta \iota(\alpha)}
\]
for $\delta \in \Dx$ such that $\delta^2 = d$ and $\varphi$ a norm
similarity of $J$ over $\D$ with multiplier $-1/\delta$.

Let $b$, $t$ and $b'$, $t'$ be triple systems associated to $\B :=
\B(J, \D\times \D)$ and $\B' := \B(J, \D\times \D, -\delta)$ given by
equations (\ref{traceform}) and (\ref{normform}).  Then define 
\[
b'_\delta :=  -\frac1{\delta} b' \hbox{\ and\ } t'_\delta := -\frac1{\delta} t'.
\]
These also define a triple system associated to $\B'$.

We have a similarity $f \!: \B \longto \B'$ as given in 
(\ref{fsim}).  Thus $\eta$ induces an $\iota$-semilinear map $\pi := f
\eta \inv{f}$ of $\B'$ given by
\[
\pi \basmat = 
\tbtmat{\iota (\beta)}{\varphi^\dagger \iota(j')}{-\varphi
  \iota(j)}{\iota (\alpha)}. 
\]
Note that $b'_\delta (\pi x, \pi y) = \iota b'_\delta (x, y)$ and similarly for
$t'_\delta$.  Thus $f$ is an isometry between the triple systems $b$, $t$
on $\B$ and $b'_\delta$, $t'_\delta$ on $\B'$ and so these triple systems
restrict to be $F$-isomorphic on $\B^\eta$ and ${\B'}^\pi$, where
$\B^\eta$ and ${\B'}^\pi$ denote the $F$-subspaces fixed by $\eta$ and
$\pi$ respectively.

A direct check (making use of Lemma \ref{jordsim}) then shows that $\pi$ is
actually an $\iota$-semilinear automorphism of $\B'$ as a Brown algebra.
Therefore, ${\B'}^\pi$ has the structure of a Brown algebra over $F$
with associated triple system $\M$.
\end{pf}


\begin{thm} \label{e7thm}
\begin{enumerate}
\item For $\B$ a Brown algebra over $F$, $\Inv(\B)$ is a
simply connected group of type $E_7$ over $F$ with trivial Tits
algebras.  Every 
simply connected group of type $E_7$ with trivial Tits algebras is obtained in
this way.
\item $\Inv(\B) \cong \Inv(\B')$ if and only if $\B \sim \B'$.
\end{enumerate}
\end{thm}

\begin{pf}
(1): Since $\Inv(\B) = \Inv(\M)$ for some nondegenerate Freudenthal
triple system 
$\M$ over $F$, it is a simply connected group of type
$E_7$ with trivial Tits algebras by \ref{FTSthm}, which finishes the
first statement.  Also by \ref{FTSthm}, if one has a simply
connected group of type $E_7$ with trivial Tits algebras, then it
is isomorphic to $\Inv(\M)$ for some 
nondegenerate Freudenthal triple system over $F$ and we are done by
Lemma \ref{goodconst}.  

(2): $\If$ is clear, so we show $\Implies$.  Brown algebras are
classified up to similarity by the 
set $H^1(F, \GInv(\B))$ and by (1) simply connected groups of type $E_7$ are
classified by $H^1(F, \PInv(\B))$.  The short exact sequence 
\[
1 \longto \GmF \longto \GInv(\B) \longto \PInv(\B) \longto 1
\]
induces an
exact sequence on cohomology
\[
H^1(F, \GmF) \longto H^1(F, \GInv(\B)) \longto H^1(F, \PInv(\B))
\]
where the last map sends $\B'$ to $\Inv(\B')$.  By
Hilbert's Theorem 90, this second map has trivial kernel.  Thus $\B
\sim \B'$.
\end{pf}


Several authors have studied Freudenthal triple systems ``in disguise''
as ternary algebras or something similar, as  
in \cite{Faulk:tern}, \cite{Faulk:E7geom}, \cite[\S 5]{FF}, \cite{FF:symp},
\cite{A:Jtern}, and \cite{Hein}.  It follows from \cite[\S 3]{FF:symp} 
that the ternary product that they are
concerned with is in fact the ternary product on a Brown algebra given 
by $(x,y,z) \mapsto \{ y, s_0 z, x \}$.
However, those authors study Lie
algebras associated to the ternary algebras which are not the Lie
algebras of the automorphism groups of the triple systems, which are
the groups that we are interested in.

\section{Singular elements in Brown algebras} \label{redsec}

For ease of notation, for an element $e$ in a Brown algebra $\B := \bbrown$,
we define a vector space endomorphism $U_e$ of $\B$ given by
\[
U_e\,x := \{ e, x, e \} \hbox{\ for all $x \in B$.}
\]

\begin{defn}
We say that 
an element $e$ in a Brown algebra $\bbrown$ is {\em singular}
if $e \ne 0$ and $U_e \, B \subseteq Fe$.
\end{defn}

In \cite[p.~196]{AF:CD} and \cite{Ferr:strict}, such elements were
called ``strictly regular''.  We are following the (shorter)
terminology from \cite{Coop:FF1} and \cite{Coop:E7} in that these elements are singular
with respect to the quartic form $q$ associated to the Brown algebra
(i.e., the radical of the quadratic form $v \mapsto q(e,e,v,v)$ is a
hyperplane of $B$ and contains $e$).

We will produce many examples of singular elements in a moment.

Following Freudenthal \cite[1.18]{Frd:E7.1} and \cite[p.~250]{SpV}, we
define a pairing $\qform{\, , \,} \!:J \times J \longto \EndF(J)$ by
\begin{neqn} \label{broket}
\qform{x, y} j := \frac12 \left( y \times (x \times j) - T(j, y) x -
  \frac13 T(x,y) j \right).
\end{neqn}

We immediately note that $\qform{x,y}^\ast = \qform{y,x}$ and that for
for $\psi \in \Inv(J)$, we have 
\begin{neqn} \label{conjverf}
\qform{ \psi(x), \psi^{\dagger}(y) } = \psi \qform{ x, y } \psi^{-1}.
\end{neqn}

\begin{lem} \label{broketlem}
Suppose that $j, j' \in J$ satisfy $j^\# = (j')^\# = 0$.  Then the
following are equivalent:
\begin{enumerate}
\item $j' \in j \times J$.
\item $\qform{j, j'} = 0$.
\item $\qform{j, j'}$ has image in $Fj$.
\end{enumerate}
\end{lem}

In the case where $j$ and $j'$ have trace zero,
the equivalence of (1) and (2) goes back to \cite[27.14]{Frd:E7.8}.

\medskip
\begin{pf}
(1) $\Implies$ (2): Suppose that $j' = j \times u$ for $u \in J$.
Then $T(j,j') = T(j\times j, u) = 0$ and $j' \times (j \times v) = T(u 
\times v, j) j = T(j', v) j$ for all $v \in J$ by \cite[(12)]{McC:FST}, so
$\qform{j, j'} = 0$.

(3) $\Implies$ (1):
We may certainly extend scalars so that we may assume that our base
field is separably closed and that we are in fact working inside the
split Albert $F$-algebra.  By (\ref{conjverf}) and Lemma
\ref{jordsim}, we may replace $j$ by
anything in its $\Inv(J^d)(F_s)$-orbit.  In particular, we may suppose 
that $j$ is the primitive idempotent $e_0 :=
\sjordmat{1}{0}{0}{0}{0}{0}$.  Then if we write $j' = \sbasjord$, we
observe that
\[
\qform{e_0, j'} \jordmat{0}{0}{0}{a'}{0}{0} =
\jordmat{0}{0}{0}{\frac23 \e_0 a'}{-a' \star c}{-b \star a'}
\]
for all $a' \in \Cd$.  Thus $\e_0 = b = c = 0$ and so $j' \in
\sjordmat{0}{F}{F}{\Cd}{0}{0}$.  This set is precisely $e_0 \times J^d$, 
as one can calculate from the explicit formula for $\times$
given in \cite[p.~358]{Jac:J}.
\end{pf}

\begin{lem} \label{sreglem}
{\rm (Cf.~\cite[6.1]{Ferr:strict})}
Let $J$ be any Albert algebra and set
$e := \sbasmat$ in $\B(J, F \times F)$.
Then $e$ is singular if and only if
\begin{enumerate}
\item $T(j, j') = 3 \alpha \beta$,
\item $(j')^{\#} = \alpha j$, 
\item $j^{\#} = \beta j'$, and
\item $\qform{j,j'} = 0$ 
\end{enumerate}

If $\alpha$ or $\beta$ is nonzero, then conditions (1) through (3)
imply (4).
\end{lem}

\begin{pf}
Direct computation shows that if (1) through (4) hold, then $e$ is
singular, so we suppose that $e$ is singular and show the converse.

If $\alpha$ or $\beta$ is nonzero, then by symmetry we may
assume that $\alpha \ne 0$.  We will deal with the $\alpha = \beta =
0$ case at the end.

Since $e$ is singular,
\begin{neqn} \label{ueozmat}
U_e \zomat = \tbtmat{2 \alpha^2}{2(j')^{\#}}{2 \alpha
j'}{T(j, j') - \alpha \beta},
\end{neqn}
 lies in $Fe$, so it must be $2 \alpha e$.  Hence conditions (1) and
 (2) hold.  Also, 
\[
U_e \ozmat = \tbtmat{2\alpha \beta}{2\beta j}{2 j^{\#}}{2 \beta^2}
\]
lies in $Fe$, so it must be $2 \beta e$.  Thus condition (3) holds.

Now we show that conditions (1) through (3) imply (4).
For any $x, y \in J$, 
\[
x^{\#} \times (x \times y) = N(x) y + T(x^{\#}, y) x
\]
by \cite[p.~496, (10)]{McC:FST}.  Since $N(j') = \frac13 T(j',
(j')^{\#}) = \alpha^2 \beta$,
\[
j \times (j' \times k') = \frac1{\alpha}(j')^{\#} \times (j'
\times k') = \alpha \beta k' + T(j, k') j'
\]
for all $k' \in J$.  Then since (3) holds, we have $0 = \qform{j', j}
= \qform{j, j'}^\ast$, hence (4).

\smallskip
Suppose now that $\alpha = \beta = 0$.   Then 
\begin{neqn} \label{qlinked} 
U_e \tbtmat{0}{k}{0}{0} = \tbtmat{0}{2 j' \times (j \times k) - k T(j, 
  j')}{\ast}{2 T(k, j^{\#})}
\end{neqn}
for all $k \in J$.  Since $T$ is nondegenerate, $j^{\#} = 0$, hence
(3).  The analogous observation with $U_e \stbtmat{0}{0}{k'}{0}$
proves (2) and considering 
(\ref{ueozmat}) demonstrates (1).  Finally, looking at equation
(\ref{qlinked}) again, we see that $\qform{j,j'}$ has image contained
in $Fj$, so by Lemma \ref{broketlem} we have (4).
\end{pf}

\begin{eg}
In $\B(J, F \times F)$, the elements $\sozmat$, $\szomat$, and
$\stbtmat{0}{j}{0}{0}$ and $\stbtmat{0}{0}{j}{0}$ for $j \in J$ such
that $j^{\#} = 0$ are all singular.
\end{eg}

\begin{eg}
Examining the proof of Lemma \ref{broketlem}, we see that
for $j' := \sjordmat{0}{0}{0}{u_1}{0}{0}$, 
the element $\stbtmat{0}{e_0}{j'}{0}$ satisfies
conditions (1) through (3) but not (4) of Lemma \ref{sreglem}.
Now, \cite[2.7]{Frd:E7.1} appears to assert that (1) through (3) imply
(4) when 
$\alpha = \beta = 0$, which would be a contradiction.  
However, Freudenthal only considered the case where $J = 
\h_3(\C, 1)$ for $\C$ a Cayley division algebra, so that $J$ has no
nonzero nilpotent elements.  For this particular algebra, conditions
(1) through (3) do imply (4), regardless of $\alpha$ and $\beta$.
\end{eg}

We have the following application.

\begin{prop} \label{type2red}
For every Albert $F$-algebra $J$ and every quadratic \'etale
$F$-algebra $\D$, $\B(J, \D) \sim \B(J, F \times F)$.  
\end{prop}

\begin{pf}
Fix some $\delta \in \Dx$ such that $\delta^2 \in \Fx$ and $\D =
F(\delta)$.  Set $\M$ to be the triple system associated to $\B(J,
\D)$ by the formulas (\ref{bform}) and (\ref{tform}) with $s_0 =
\delta \stbtmat{1}{0}{0}{-1}$.  By the lemma,
\[
f_1 := \stbtmat{1}{1}{1}{1} \hbox{\ and\ } f_2 := \delta
\stbtmat{-1}{-1}{1}{1}
\]
are singular elements in $\M$.

The map $h \!: \M \otimes_F \D \longto \M(J) \otimes_F \D$ given by
\[
h \basmat := \tbtmat{\alpha / \delta}{\delta j}{j'}{\delta^2 \beta}
\]
is an isomorphism of triple systems.
As described in \cite[p.~95]{Brown:E7}, for $k \in J$ we have
automorphisms of $\M(J)$ given by 
\begin{neqn} \label{phidef}
\varphi_k \basmat = \tbtmat{\alpha + \beta N(k) + T(j',
  k) + T(j, k^\#) }{j + \beta k}{j' + j \times k + \beta k^\#}{\beta}
\end{neqn}
and
\begin{neqn} \label{psidef}
\psi_k \basmat = \tbtmat{\alpha}{j + j' \times k + \alpha k^\#}{j' +
  \alpha k}{\beta + \alpha N(k) + T(j, k) + T(j', k^\#)},
\end{neqn}
where $\varphi_k^{-1} = \varphi_{-k}$ and $\psi_k^{-1} = \psi_{-k}$.
Then
\[
m := \varphi_{-\frac1{2\delta}}\, \psi_\delta\, h \!: \M \otimes_F \D \longto 
\M(J) \otimes_F \D
\]
is an isomorphism such that 
\[
m(f_1) = \tbtmat{0}{0}{0}{8 \delta^2} \hbox{\ and\ } m(f_2) =
\tbtmat{-1}{0}{0}{0}.
\]

Now for $\iota$ the nontrivial $F$-automorphism of $\D$, consider the
1-cycle $z \in Z^1(\D/F, \Inv(\M(J)))$ given by $z_\iota := m (\varpi
\otimes \iota) m^{-1} (1 \otimes \iota)^{-1}$.  It fixes the diagonal
matrices in $\M(J)$ elementwise, so by \cite[7.5]{Ferr:strict} it is
equal to $f_\phi$ for some $\phi \in \Inv(J)(\D)$, in the notation of
(\ref{fdef}).  Then the obvious computation shows that $\phi$ is the 
identity, so that $z$ is the trivial cocycle and $\M \cong \M(J)$. 
\end{pf}

\begin{borel*} \label{nonredrmk} 
Not all Brown algebras contain singular elements.  Having a singular
element  corresponds to 
one (equivalently, all) of the triple systems associated to it being
of the form $\M(J)$ for some $J$.  Ferrar gives an example
in \cite[p.~330]{Ferr:strict} and \cite[pp.~64, 65]{Ferr:E6} of a
nondegenerate Freudenthal triple system $\M$ over a  
field of transcendence degree 4 over $\Q$ such that $\M$ is not of
such a form. By
Lemma \ref{goodconst}, $\M$ arises from some Brown algebra $\B$ which
then contains no singular elements, so in particular is of type 2 and
not of the form $\B(J, \D)$ for any $J$ or $\D$ by \ref{type2red}.
\end{borel*}

But there is something we can say in a special case.  The following
result was suggested by Markus Rost: 
\begin{prop}
Suppose that $\B$ is a Brown $F$-algebra of type 2 with inner
extension $\D$ such that $\B \otimes_F \D$ is split.  Then $\B$
contains a singular element.
\end{prop}

\begin{pf}
Let $B^q := \B(J^d, \D)$ be the quasi-split Brown algebra with inner
extension $\D$.  Then there is some class $(f) \in H^1(\D/F,
\autp(\B^q))$ which corresponds to $\B$, which must be of the form
\[
f_\iota \basmat = \tbtmat{\alpha}{\varphi(j)}{\varphi^\dag(j')}{\beta}
\]
where $\iota$ is the nontrivial $F$-automorphism of $\D$ and $\varphi
\in \Inv(J^d)(\D)$.  Since $f$ is a 1-cocycle, $\varphi \iota
\varphi^\dagger \iota$ is the identity in $\Inv(J^d)(\D)$.  This is
the situation addressed in \cite[p.~65, Lem.~3]{Ferr:E6}, and the
proof of that result shows that for any element $u \in J^d$ such that
$u^\# = 0$, we can modify $f$ and so assume that $\varphi \iota(u) =
\delta u$ for some $\delta \in \Dx$.  We pick such a $u$ so that
$T(u,u) = 0$ and $\qform{u,u} = 0$, for example,
\[
u := \sjordmat{0}{0}{0}{u_1}{0}{0}.
\]
Then the element
\[
\stbtmat{0}{u}{\iota(\delta)^{-1}u}{0}
\]
is fixed by $f_\iota \iota$ and so is an element of $\B$.  it is
singular by Lemma \ref{sreglem}.
\end{pf}

\section{Inner ideals} \label{idealsec}

\begin{defn} 
A vector subspace $I$ of a structurable algebra $(A, {-})$ is said to
be an {\em inner ideal} if $U_e (A) \subseteq I$ for all $e \in I$.
We say that $I$ is {\em proper} if $I \ne A$.
\end{defn}

\begin{eg}[McCrimmon]
One says that an element $d$ in an Albert $F$-algebra $J$ is {\em of
  rank one} if $d^{\#} = 0$.  We say that an $F$-subspace $V$ of $J$
is {\em totally singular} if every element of $V$ is of rank one.  By
\cite[p.~467, Thm.~8]{McC:inn}, the proper inner ideals in an Albert
algebra are the totally singular subspaces and the subspaces of the
form $d \times J$ for some $d$ of rank one.
We will call this last sort of subspace a {\em hyperline}, following
Tits' terminology from \cite{Ti:Rsp}.
\end{eg}

\begin{eg} \label{singeg}
For $V$ a totally singular subspace as in the preceding example,
$\stbtmat{F}{0}{V}{0}$ is an inner ideal in $\B(J, F \times F)$ since
every element is singular by Lemma \ref{sreglem}.  
\end{eg}

We say that an inner ideal in a Brown algebra is {\em singular} if it
consists of  
singular elements.  We call such inner ideals {\em singular ideals}
for short.  (Of course, any subspace of a Brown algebra consisting
of singular elements is automatically a singular ideal.)

\begin{eg} \label{nsingeg}
Suppose that $d \in J$ satisfies $d^{\#} = 0$.  Then 
\[
I := \tbtmat{F}{Fd}{d \times J}{0}
\]
is an inner ideal of $\B(J, F \times F)$.  Moreover, it is not
singular, since the element $\stbtmat{1}{d}{0}{0}$ 
is in $I$ and it is not singular by Lemma \ref{sreglem}.
(After we have proven Lemma
\ref{hypsize}, we will also know that $I$ is 12-dimensional.)
\end{eg}

We want to classify the inner ideals in a Brown algebra $\B$ enough to
describe the homogeneous projective varieties associated to groups of
type $E_7$ in terms of such ideals.  Our classification is going to
rest upon first understanding the inner ideals in Albert algebras a
bit better.

\begin{lem}\label{sharplem}
The set $\{ j^{\#} \mid j \in J^d \}$ spans $J^d$ over any separably
closed field.
\end{lem}

\begin{pf}
In the notation of (\ref{rednot}), if one sets precisely one entry to
be nonzero (and the entry to be $u_i$ for some $i$ if it is $a$, $b$,
or $c$) then one gets a rank one element and
these elements span $J^d$. 
Thus it suffices to prove
that for every $x$ of rank one, there is some $j \in J^d$ with $j^{\#}
\in Fx$.  As discussed in \cite[p.~70]{Jac:J3}, for $x \ne 0$ of rank one, $x$
is either (1) a primitive idempotent or it satisfies (2) $T(x,1) = 0$
and $x^2 = 0$.

(1):  Let $e_0$, $e_1$, and $e_2$ be a triple of
orthogonal idempotents in $J^d$ such that $1 = e_0 + e_1 + e_2$. The
group $\aut(J^d)(F)$ acts transitively on the primitive 
idempotents of $J^d$ by the Coordinatization Theorem
\cite[p.~137]{Jac:J}, so we may assume that $x = e_0$.  Then $(e_1 +
e_2)^{\#} = e_0$.

(2): The group $\aut(J^d)(F_s)$ acts transitively on such rank one
elements by \cite[28.22]{Frd:E7.8}, so we may assume that
\[
x = \jordmat{0}{0}{0}{0}{u_1}{0} \hbox{, and so\ } j :=
\jordmat{0}{0}{0}{u_1}{0}{-u_4} 
\]
satisfies $j^{\#} = x$.
\end{pf}

\begin{lem} \label{hypsize}
Any hyperline in an Albert algebra is 10-dimensional.
\end{lem}

\begin{pf}
We may certainly extend scalars to assume that our base field is
separably closed and that the Albert algebra is quasi-split.  Suppose
our hyperline is $d \times J^d$.  Then
there is some norm isometry $\varphi$ such that $\varphi(d)$ is the 
primitive idempotent
$e_0 := \sjordmat{1}{0}{0}{0}{0}{0}$.  Of course,
\[
\dim (x \times J) = \dim \varphi^{\dagger}(x \times J) = \dim (e_0
\times J),
\]
But we already know what $e_0 \times J$ is from the proof of Lemma
\ref{broketlem}, and that it is 10-dimensional.
\end{pf}

What follows is our last preparatory lemma concerning Albert algebras.
 It may seem mysterious now, but after Lemma \ref{e6flag} we
 will see that it is precisely a simple algebraic interpretation of
 the nontrivial automorphism on the Dynkin diagram of type $E_6$, or,
 if you prefer, of the natural duality in the spherical building
 of type $E_6$
 a.k.a.~a Hjelmslev-Moufang plane.  It is stronger than what we
 actually need for the rest of this section, but it will all be of use later.

\begin{dualitylem}\label{dualitylem}
Let $J$ be an Albert $F$-algebra.  The map on subspaces of $J$ which takes
a subspace $W$ of $J$ to
\[
\{ j \in J \mid \hbox{$\qform{ w, j } = 0$ for all $w \in W$} \}
\]
induces one-to-one correspondences
\[
\begin{array}{rcl}
\hbox{3-dim'l t.~singular subspaces}&\leftrightarrow&\hbox{3-dim'l t.~singular subspaces}\\
\hbox{2-dim'l t.~singular subspaces}&\leftrightarrow&\hbox{5-dim'l
  {\em maximal} t.~singular subspaces} \\
\hbox{1-dim'l t.~singular subspaces}&\leftrightarrow&\hbox{hyperlines}
\end{array}
\]
The bottom correspondence is given by 
\[
Fd \leftrightarrow d \times J.
\]
\end{dualitylem}

In particular this gives us the simple fact that if $d \times J = d'
\times J$ for rank one elements $d$ and $d'$, then $Fd = Fd'$.

\medskip
\begin{pf}
We may certainly assume that our base field $F$ is separably closed and
that our Albert algebra is split.  Then by \cite[3.2, 3.12, 3.14]{SpV}
the group of norm isometries acts transitively on each of the six 
kinds of subspaces specified.  So all we really need to do is produce
an example of a pair of subspaces $(W,W')$ which are sent to each
other by the specified map.  For if $V$ is another subspace of the
same kind as $W$, there is some norm isometry $\psi$ with
$\psi(V) = W$.  By (\ref{conjverf}), our correspondence map will
then send $V$ to $\psi^\ast(W')$, which in turn is itself sent to
$V$, proving that the correspondence is indeed a bijection.

Consider the primitive idempotent $e_0$ from the proof of Lemma
\ref{hypsize}.  We will show that $(Fe_0, e_0 \times J)$ provides a
pair which is an example of the last correspondence.  By Lemma
\ref{broketlem}, $Fe_0$ is sent to $e_0 \times J$, so we must check
the other direction.  We can define other primitive idempotents $e_1$ and
$e_2$ such that $e_i$ has all entries zero except for a one in the
$(i+1, i+1)$-position.  
Since $e_1, e_2 \in e_0 \times J$, by Lemma \ref{broketlem} $e_0
\times J$ is sent to a subspace $W$ 
which is 
contained in
\[
(e_1 \times J) \cap (e_2 \times J) = \sjordmat{F}{0}{F}{0}{\Cd}{0} \cap
\sjordmat {F}{F}{0}{0}{0}{\Cd} = F e_0.
\]
But $W$ certainly contains $F e_0$, so $e_0 \times J$ is sent to $W =
Fe_0$.

Note that we have now even proven the claim about the form of the
last correspondence, since for a norm isometry $\psi$ such that $F
\psi(d) = F e_0$, $Fd$ is sent to $\psi^\ast(e_0 \times J) =
\psi^{-1}(e_0) \times J = d \times J$.

Now set 
\[
W := \jordmat{0}{0}{F}{Fu_1}{0}{0} \hbox{\ and\ } W' :=
\jordmat{F}{0}{0}{0}{0}{\C \star u_1}.
\]
It is quickly checked that $W$ is a 2-dimensional totally singular
subspace and that $W'$ is a 5-dimensional maximal totally singular
subspace.
Similarly to what we just did for the last correspondence, we can
apply Lemma \ref{broketlem} to rewrite the correspondence map as
\[
W \mapsto \cap_{w \in W} w \times J.
\]
Then it is just a tedious computation to see that the pair $(W,W')$ we 
have just defined are indeed sent to each other.

Finally,
\begin{neqn} \label{tredim}
W := W' := \jordmat{0}{0}{0}{Fu_1}{Fu_2}{Fu_5}
\end{neqn}
provide a pair $(W,W')$ of 3-dimensional totally singular subspaces
which are mapped to each other by the correspondence map.  (There are
other choices for $(W,W')$ which would be easier to check here, but we 
will use this subspace again at the end of Section \ref{flagsec}.  We
also caution
the reader that this pair is atypical in that a 3-dimensional totally
singular space is not necessarily sent to itself by our
correspondence map.) 
\end{pf}

\begin{rmk}
Incidentally, the other sort of maximal totally singular subspaces are
6-dimensional and $\Inv(J^d)$ acts transitively on the set of such
subspaces \cite[3.14]{SpV}.  A particular example is given by
\[
\jordmat{0}{0}{F}{\Cd \star u_1}{F u_1}{0}.
\]
Since this subspace is sent to zero by the correspondence map above,
so are all 6-dimensional totally singular subspaces.
\end{rmk}

\begin{lem} \label{bigideallem}
Let $I$ be an inner ideal in $\B(J, F \times F)$.  If $I$ contains
$\stbtmat{0}{0}{J}{0}$ or $\stbtmat{0}{J}{0}{0}$, then $I$ is all of $\B$.
\end{lem}

\begin{pf}
Since the property of being proper is invariant under scalar
extension, we may assume that our base field is separably closed and
$J$ is split (i.e., $J = J^d$).

In the first case, $I$ contains
\[
U_{\stbtmat{0}{0}{j'}{0}} \zomat = \tbtmat{0}{2(j')^{\#}}{0}{0}
\]
for all $j' \in J^d$.  Since $(J^d)^{\#}$ spans $J^d$ by Lemma \ref{sharplem}, $I$ contains
$\stbtmat{0}{J^d}{J^d}{0}$.  

Then $I$ also contains
\[
U_{\stbtmat{0}{j}{j'}{0}} \zomat = \tbtmat{0}{2(j')^{\#}}{0}{T(j,j')},
\]
so $I$ contains $\stbtmat{0}{J^d}{J^d}{F}$.  Then
\[
U_{\stbtmat{0}{j}{j'}{0}} \ozmat = \tbtmat{T(j,j')}{0}{2j^{\#}}{0} \in I,
\]
so we are done with the first case.

The proof in the second case is symmetric to the first case.
\end{pf}

We also observe in the following proposition that whether or not a subspace
is a singular or an 
inner ideal is already determined by its characteristics in an
associated triple system, so that $\Inv(\B)$ takes inner
(resp.~singular) ideals to inner (resp.~singular) ideals.

\begin{prop} \label{e7obs} \label{e7prop}
Let $\B = \bbrown$ be a Brown algebra with $b$, $t$ an associated
triple system.  Then a vector subspace $I$ is an inner ideal of $\B$ if and 
only if 
\[
t(I, I, B) \subseteq I.
\]
It is singular if and only if
\[
t(u, v, z) = b(z,v) u + b(z,u) v \hbox{\ for all $u, v \in I$ and $z \in B$}.
\]
\end{prop}

\begin{pf}
Observe that 
\[
t(u,v,z) = \frac12 \left( t(u,z,v) + t(v,z,u) \right) = \{ u, s_0 z, v
\} + \{ v, s_0 z, u \} - b(z,v) u - b(z,u) v.
\]
$I$ is an inner ideal if and only if the sum of the two brace terms on
the left-hand side is in $I$ for all $u, v \in I$, so this proves the
first equivalence.

By \cite[p.~317]{Ferr:strict}, an element $u \in B$ is strictly
regular if and only if $t(u,u,y) = 2 b(y, u) u$ for all $y \in B$.  Since
\[
t(u + v, u + v, z) = t(u, u, z) + t(v, v, z) + 2 t(u, v, z),
\]
if $u$, $v$, and $u + v$ are strictly regular, then 
\begin{neqn} \label{sreq}
t(u, v, z) = b(y, u + v) (u + v) - b(y,u)u - b(y,v)v =
b(y, u) v + b(y, v) u.
\end{neqn}
for all $z \in B$.  Conversely, if (\ref{sreq}) holds for all $u, v
\in I$, then every element of $I$ is singular.
\end{pf}

\begin{thm} \label{transthm}
If a proper inner ideal of a Brown algebra $\B$ contains a singular
element, then it is in the $\Inv(\B)(F_s)$-orbit of an inner ideal as in  
Example \ref{singeg} or Example \ref{nsingeg}.  In particular, it has
dimension $\le 7$ and is singular or is $12$-dimensional and is not singular.
\end{thm}

\begin{pf}
We may assume that our base field $F$ is separably closed and so that
our algebra is split.
Since $\Inv(\B)(F_s)$ acts transitively on singular elements
\cite[7.7]{Ferr:strict}, we may assume that the ideal $I$ contains $e
:= \sozmat$.

Since $e \in I$, we may extend $e$ to a basis $x_2, \ldots,
x_n$ of $I$ such that 
\[
x_i := \tbtmat{0}{j_i}{j'_i}{\beta_i}
\]
for $2 \le i \le n = \dim_F I$.

\newcommand{\yto}{\tbtmat{0}{0}{k'}{0}}

First, consider
\[
\begin{array}{rcl}
\left( U_{e + x_i} - U_e - U_{x_i} \right) \yto &=& \left\{ e, \yto, x_i \right\} + \left\{ x_i, \yto, e \right\} \\
&=& \tbtmat{0}{0}{\beta_i k'}{0} \in I.
\end{array}
\]
Thus, if any $\beta_i$ is nonzero, the ideal would have to contain
$\stbtmat{F}{0}{J}{0}$, and so by Lemma \ref{bigideallem} it would not be
proper.  Thus $\beta_i = 0$ for all $i$.

If $I$ is singular, then take $x := \stbtmat{0}{j}{j'}{0}$ in $I$.
Since $x + \alpha e \in I$ for all $\alpha \in F$, it is singular.
By Lemma
\ref{sreglem}, $(j')^{\#} = \alpha j$ for all $\alpha \in \Fx$, so $j
= 0$.  Then $I$ is one of the ideals described in Example
\ref{singeg}.

Otherwise, 
\[
\left( U_{x_i + x_\ell} - U_{x_i} - U_{x_\ell} \right) \tbtmat{0}{k}{0}{0} = 
\tbtmat{0}{\ast}{\ast}{2 T(k, (j_i
\times j_\ell))} \in I.
\]
Since $T$ is nondegenerate, $j_i \times j_\ell = 0$ for all $2 \le i,
\ell \le n$.  So if we write $I \subseteq \stbtmat{F}{W}{V}{0}$ where
$W$ and $V$ are the projections of $I$ into $J$ on the off-diagonal
entries, we have shown 
that $w^{\#} = 0$ for all $w \in W$.  

Let $x := \stbtmat{\alpha}{w}{v}{0}$ represent an arbitrary element of
$I$.  Then
\[
U_{x - \alpha e} \zomat = \tbtmat{0}{2 v^{\#}}{0}{T(v, w)} \in I,
\]
so $T(v,w) = 0$ and $v^{\#} \in W$.  Also,
\begin{neqn} \label{ualpha}
U_x \tbtmat{0}{k}{0}{0} = \tbtmat{\ast}{\ast}{2 \alpha w \times k - 2
v^{\#} \times k + 2 v T(v, k)}{0}.
\end{neqn}
Setting $\alpha = 0$, we see that $V^{\#} \times J \subseteq V$.  (We
will use the $\alpha \ne 0$ case in a moment.)  Thus $V$ is an
inner ideal of $J$ by \cite[p.~467]{McC:inn}.  We consider the
following cases: (1) $V = J$, (2) $V$ is a hyperline, or (3) $V^{\#} = 0$.

(1) cannot occur because in that case $V^{\#} = J$ by Lemma
\ref{sharplem}, and we already know that $V^{\#} \subseteq W$ and
$W^{\#} = 0$.

In case (2), $V = d \times J$ for some $d$ of rank one.  However,
letting $\alpha$ be arbitrary  
in (\ref{ualpha}) and keeping in mind that $V^\# \times J \subseteq
V$, we see that in fact $W \times J \subseteq V = d \times J$.
Since $W^\# = 0$, $w \times J = w' \times J$ if and only if $w \in
Fw'$ by the Duality Lemma \ref{dualitylem}.
The fact that hyperlines are maximal proper inner ideals of $J$
\cite[p.~467, Thm.~8]{McC:inn} tells us that $W \subseteq Fd$.
Since $Fd = V^{\#} \subseteq W$, the ideal is as in
Example \ref{nsingeg}.

Finally, we examine case (3).  There, since $V^{\#} =
0$, $V$ is at most 6-dimensional by \cite[3.14]{SpV}.  However, $W^\#
= 0$ and $W \times J \subseteq V$.  Since hyperlines are
10-dimensional, $W = 0$.
Thus $I$ is as in Example \ref{singeg}.
\end{pf}

\begin{thm} \label{bigidealthm}
Any proper inner ideal in a Brown algebra has dimension at most $12$.  If
it is $12$-dimensional, then it is in the $\Inv(\B)(F_s)$-orbit of an
inner ideal as in Example \ref{nsingeg}.
\end{thm}

\begin{pf}
Clearly we may assume that our base field is separably closed and so that
$\B = \B(J, F \times F)$ for $F = F_s$.  Suppose that $\dim_F I \ge
12$ and pick a
basis $x_1$, $x_2$, $\ldots$, $x_n$ of $I$ such
that 
\[
x_i := \tbtmat{\alpha_i}{j_i}{j'_i}{\beta_i}
\]
where $\beta_i = 0$ for $2 \le i \le n$.  Set $c := 1$
if $\beta_1 = 
0$, otherwise set $c := 2$.  Finally, set $I'$ to be the span of $x_i$
for $c \le i \le n$.  Let $W$ denote the span of $j_i$ for $c \le i
\le n$.  

For $x := \stbtmat{\alpha}{j}{j'}{0}$ an arbitrary element of $I'$,
\[
U_x \tbtmat{0}{0}{k'}{0} = \tbtmat{\ast}{2 j T(j, k') - 2 k' \times
j^{\#}}{\ast}{0} \in I.
\]
In fact, it lies in $I'$.  Thus $W^{\#} \times J \subseteq W$ and $W$ is
an inner ideal of $J$, so (1) $W = J$, (2) $W$ is a hyperline, 
or (3) $W^{\#} = 0$.

In case (1),
\[
U_x \ozmat = \tbtmat{0}{0}{2 j^{\#}}{0} \in I',
\]
so by Lemma \ref{sharplem} $I$ contains $\stbtmat{0}{0}{J}{0}$.  By
Lemma \ref{bigideallem}, $I$ is not proper, so this is a
contradiction.

In case (2) $W$ is 10-dimensional and in case (3) $W$ is at most
6-dimensional by \cite[3.14]{SpV}.  Thus we may rewrite the basis so
that $j_i = 0$ for $i \ge c + 10$ (which is $\le 12$).  Consider 
\[
x := x_{c+10} = \tbtmat{\alpha}{0}{v}{0},
\]
for $v := j'_{c+10}$. If $v = 0$, then $x \ne 0$ would imply that $I$
contains a singular element and so we would be done by
Theorem \ref{transthm}.  So we may assume that $v \ne 0$.
If $v^{\#} = 0$ then
$x$ is singular and we are done.
Otherwise, consider
\[
x' := \frac12 U_x \zomat - \alpha x = \tbtmat{0}{v^{\#}}{0}{0} \in I.
\]
If $v^{\#\#} = 0$ then $x'$ is singular and we are again done, so we
may assume that 
$N(v) v = v^{\#\#} \ne 0$.  Then
\[
x''_k := \frac12 U_{x'} \tbtmat{0}{k}{0}{0} = \tbtmat{0}{0}{0}{T(k,
  v)N(v)} \in I 
\]
for all $k \in J$.  Since $T$ is nondegenerate, $x''_k \in I$ is
nonzero and singular for some choice of $k$.
\end{pf}

\section{Flag varieties for groups of type $E_7$} \label{flagsec}
\setlength{\unitlength}{1cm}
\newcommand{\darkradE}{0.115}

In Section \ref{idealsec} I promised that the homogeneous projective
varieties associated to the 
group $\Inv(\B)$ of type $E_7$ can be described in terms of the inner ideals
of $\B$.  We will fulfill this promise in \ref{e7flag}, but first we
must  set up some notation. 

\begin{borel}{Background on flag varieties}\label{flagback}

Given a split maximal torus $T$ in a simple affine algebraic group $G$, we
fix a set of simple roots $\Delta = \{ \alpha_1, \ldots, \alpha_r \}$
of $G$ with respect to $T$.  For each $\alpha_i$, there is a uniquely
defined root group $U_{\alpha_i}$ lying in $G$ \cite[13.18]{Borel}.
For any subset $\Theta = \{ i_1, i_2, \ldots, i_n \}$ of $\Delta$, we
define a parabolic subgroup of $G$ by
\[
P(\Theta) := \qform{T, \{ U_\alpha \mid \alpha \in \D \}, \{
  U_{-\alpha} \mid \alpha \not\in \D \} }
\]
and an associated flag variety by
\[
X(\Theta) := G/P(\Theta).
\]
Thus $P(\Delta)$ is a Borel subgroup of $G$ and $P(\emptyset) = G$.
This notation is similar to \cite{MPW1} and \cite{G:flag} and is
opposite that of \cite[4.2]{BoTi} and \cite{KR}.

We are interested in the group $G = \Inv(A)$ for $A$ an Albert or
Brown $F$-algebra.  Each variety $X(i)$ is going to have $F$-points
corresponding to certain subspaces of $A$ which we will call {\em
  $i$-spaces}.  We will also define a symmetric and reflexive binary
relation called {\em incidence} between $i$-spaces and $j$-spaces for
all $i$ and $j$.  Two $i$-spaces will be incident if and only if they
are the same.  Then the other flag varieties will be of the form
$X(\Theta)$ and have $F$-points
\[
\left\{ (V_1, \ldots, V_n) \in X(i_1)(F) \times \cdots \times X(i_n)(F)
    \>\left|\>
\parbox{3.5cm}{$V_i$ is incident to $V_j$ for all $1 \le i, j \le n$}
\right. \right\}.
\]
\end{borel}

\subsection*{Flag varieties for $E_6$}
As a prelude to describing the flag varieties for $\Inv(\B)$, we
recall the description of the flag varieties of $\Inv(J)$ for $J$ an
Albert algebra.  We number the Dynkin diagram of $\Inv(J)$ as
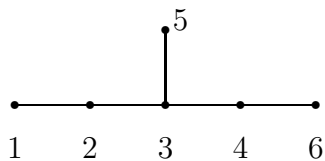
\[
\begin{picture}(7,3)
    \multiput(1,1)(1,0){5}{\circle*{\darkradE}}
    \put(3,2){\circle*{\darkradE}}

    \put(1,1){\line(1,0){4}}
    \put(3,2){\line(0,-1){1}}
    
    \put(1,0.3){\makebox(0,0.4)[b]{$1$}}
    \put(2,0.3){\makebox(0,0.4)[b]{$2$}}
    \put(3,0.3){\makebox(0,0.4)[b]{$3$}}
    \put(4,0.3){\makebox(0,0.4)[b]{$4$}}
    \put(5,0.3){\makebox(0,0.4)[b]{$6$}}
    \put(3.2,2){\makebox(0,0.4)[b]{$5$}}

\end{picture}
\]

\begin{thm} \label{e6flag}
Let $J$ be an Albert $F$-algebra.  The $i$-spaces for $\Inv(J)$ are
the $i$-dimensional totally singular subspaces of $J$ for $i = 1, 2,
3$.  The $4$-spaces are the $5$-dimensional {\em maximal} totally
singular subspaces of $J$.  The $5$-spaces are the $6$-dimensional
totally singular subspaces, and the $6$-spaces are the subspaces of
the form $d \times J$ for $d \in J$ of rank one.

Incidence is defined by inclusion except for the following: A
$4$-space and a $5$-space are incident if and only if their
intersection is $3$-dimensional.  A $6$-space and a $5$-space are
incident if and only if their intersection is $5$-dimensional (in
which case it is necessarily a nonmaximal totally singular subspace of
$J$).
\end{thm}

\begin{pf}
For the purposes of the proof we define an $i$-space to be a subspace
of $J$ as specified in the theorem statement, and we will show that one
can identify $X(i)(F)$ with the set of $i$-spaces.  

Let $X_i$ denote
the functor mapping field extensions of $F$ to sets of subspaces of
$J$ such that $X_i(K)$ is the set of $i$-spaces of $J \otimes_F K$.
Then certainly $X_i$ is a projective variety for $i = 1, 2, 3, 5$.
Duality for points and hyperlines (\ref{dualitylem}) shows that $X_6$ and $X_4$
are projective varieties since $X_1$ and $X_2$ are.

The transitivity of the natural $\Inv(\B)$-action on $X_i$ is given by
\cite[3.2]{SpV} for $i = 1$, by \cite[3.12]{SpV} for $i = 2, 3$, and
by \cite[3.14]{SpV} for $i = 4, 5$.  The transitivity of the action on
$X_6$ follows by the definition of a hyperline and the transitivity of
the action on $X_1$.
The incidence relations and the
associations with simple roots are in \cite[3.2]{Ti:Rsp}, but we will
produce a set of simple roots explicitly because we will need them
later when we consider $E_7$.  All of the material we develop here
will also see use in Section \ref{cisec}. 

In order to describe the associations with simple roots, we extend
scalars to split $J$.  

Consider the split Cayley algebra $\Cd$ with
hyperbolic norm form $\n$ and multiplication $\star$ as
described in the introduction.  We will define two algebraic groups
$\Rel\Cdn$ and $\Spin\Cdn$ associated to $\Cdn$.  First, consider the
(connected, reductive) algebraic group $\GOP\Cdn$ whose $F$-points are
\[
\GOP\Cdn(F) := \{ f \in \EndF(\Cd) \mid \hbox{$\s(f)f \in \Fx$ and $(\det
  f)^4 = \s(f)f$} \},
\]
where $\s$ is the involution on $\EndF(\Cd)$ which is adjoint to $\n$.
For $f \in\GOP\Cdn(F)$, we write $\mu(f) := \s(f)f$, the {\em multiplier}
of $f$.  We define $\Rel\Cdn$ to be the algebraic group whose
$F$-points are the {\em related triples} in $\trio{\GOP\Cdn}$, i.e., those
triples $\ut := \trion{t}$ such that 
\[
\mu(t_i)^{-1} t_i(x \star y) = t_{i+2}(x) \star t_{i+1}(y)
\]
for $i = 0,1,2$ (subscripts taken modulo 3) and all $x, y \in \Cd$.
We write $\OP\Cdn$ for the closed subgroup of $\GOP\Cdn$ consisting of
those $f$'s with $\det f = 1$ and define $\Spin\Cdn$ to be the
algebraic group consisting of the related triples in $\trio{\OP\Cdn}$.
Now $\Rel\Cdn$ injects into $\Inv(J^d)$ by sending $\ut$ to $g_\ut$,
which is defined by
\begin{neqn} \label{gdef}
g_\ut \basjord := \jordmat{\mu(t_0)^{-1} \e_0}{\mu(t_1)^{-1}\e_1}{\mu(t_2)^{-1}\e_2}{t_0 a}{t_1 b}{t_2 c}.
\end{neqn}
This map restricts to an injection of $\Spin\Cdn$ into $\aut(J^d)$.
The subset of $\Spin\Cdn$ consisting of triples $\ut$ such that $t_i$
is diagonal for all $i$ forms a rank 4 split torus in $\Spin\Cdn$
\cite[1.6]{G:iso}, and the image of this torus under the composition
$\Spin\Cdn \injects \aut(J^d) \injects \Inv(J^d)$ provides a rank 4 split
torus $S_4$ in $\Inv(J^d)$.

For $\ul = \trion{\lambda} \in
\trio{(\Fx)}$ such that $\lambda_0 \lambda_1 \lambda_2 = 1$, we have a
map $S_{\ul} \in \Inv(J^d)$ given by
\[
S_{\ul} \basjord = \jordmat{\lambda_0^{-2} \e_0}{\lambda_1^{-2}
  \e_1}{\lambda_2^{-2} e_2}{\lambda_0 a}{\lambda_1 b}{\lambda_2 c}.
\]
Let $S_2$ denote the rank 2 split torus in $\Inv(J^d)$ generated by such
maps.  Then $S_6 := S_2 \times S_4$ is a rank 6 split torus in
$\Inv(J^d)$.  We have characters $\chi_{i,j}$ defined for $0 \le i \le
2$, $1 \le j \le 8$ by setting $\chi_{i,j}$ to be trivial on $S_2$ and
to take the value of the $(j,j)$-entry of $t_i$ on $\trion{t}$.
Define $\rho_i$ to be the character which is trivial on $S_4$ and such
that
$\rho_i (S_{\trion{\lambda}}) = \lambda_i$.

A set of simple roots is given by the following, where we have written
$\omega_j := \chi_{0,j}$ for short:
\begin{ngather}
\begin{gather}
\alpha_1 = -\omega_1 - (\rho_2 - \rho_1), \label{e6roots}\\
\alpha_2 = \omega_1 - \omega_2,\quad \alpha_3 = \omega_2 - \omega_3,\quad
\alpha_4 = \omega_3 + \omega_4,\quad \alpha_5 = \omega_3 - \omega_4,
\notag \\
\alpha_6 = \chi_{2,8} - (\rho_1 - \rho_0) = -\frac12 (\omega_1 +
\omega_2 + \omega_3 + \omega_4) - (\rho_1 - \rho_0)  \notag
\end{gather}
\end{ngather}
The root groups corresponding to the root subsystem of type $D_4$
spanned by $\alpha_2$ through $\alpha_5$ are given explicitly in
\cite[4.3]{G:iso}.
The 1-dimensional root Lie algebras corresponding to the roots
$\alpha_1$ and $\alpha_6$ are $S_{(Fu_8)_{23}}$ and $S_{(Fu_8)_{12}}$ 
respectively, in the notation of \cite[p.~35]{Jac:ex}.

Consider the following spaces:
\begin{gather*}
V_1 := \sjordmat{0}{0}{F}{0}{0}{0}, \quad
V_2 := \sjordmat{0}{0}{F}{Fu_1}{0}{0},\quad
V_3 := \sjordmat{0}{0}{F}{Fu_1 + Fu_2}{0}{0}, \\
V_4 := \sjordmat{0}{0}{F}{u_1 \star \Cd}{0}{0}, \quad
V_5 := \sjordmat{0}{0}{F}{\Cd \star u_1}{Fu_1}{0}, \quad \hbox{\ and\  } \quad
V_6 := e_0 \times J = \sjordmat{0}{F}{F}{\Cd}{0}{0},
\end{gather*}
for $e_0$ the primitive idempotent as in the proof of \ref{broketlem}.
Then $V_i$ is the sort of space which is claimed to be an $i$-space.

We leave it to the reader
to verify that the stabilizer of $V_i$ in $\Inv(J^d)$ is precisely the
parabolic $P(i)$, so that in fact $X_i = X(i)$.  The only data the
reader is potentially missing is that the root Lie algebras
corresponding to the roots $-\alpha_1$ and $-\alpha_6$ are
$S_{(Fu_1)_{32}}$ and $S_{(Fu_1)_{21}}$. 
\end{pf}

\subsection*{Flag varieties for $E_7$}

We label the Dynkin diagram for $E_7$ as follows
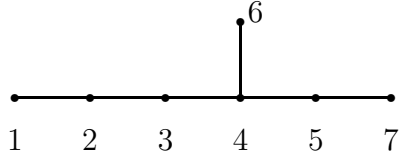
\[
\begin{picture}(7,3)
    \multiput(1,1)(1,0){6}{\circle*{\darkradE}}
    \put(4,2){\circle*{\darkradE}}

    \put(1,1){\line(1,0){5}}
    \put(4,2){\line(0,-1){1}}
    
    \put(1,0.3){\makebox(0,0.4)[b]{$1$}}
    \put(2,0.3){\makebox(0,0.4)[b]{$2$}}
    \put(3,0.3){\makebox(0,0.4)[b]{$3$}}
    \put(4,0.3){\makebox(0,0.4)[b]{$4$}}
    \put(5,0.3){\makebox(0,0.4)[b]{$5$}}
    \put(6,0.3){\makebox(0,0.4)[b]{$7$}}
    \put(4.2,2){\makebox(0,0.4)[b]{$6$}}

\end{picture}
\]
The key idea here is that nodes 2 through 7 span a Dynkin diagram of
type $E_6$, and that this corresponds to the subgroup $\Inv(J)$ of
$\Inv(\B(J))$.

\begin{thm} \label{e7flag}
Let $\B$ be a Brown $F$-algebra.  The $i$-spaces for $\Inv(\B)$ are
the $i$-dimensional singular ideals for $i = 1$,$2$,$3$,$4$.  The
$5$-spaces are the $6$-dimensional {\em maximal} singular ideals.  The
$6$-spaces are the $7$-dimensional singular ideals, and the $7$-spaces
are the $12$-dimensional inner ideals.  

Incidence is defined by inclusion except for the following: A
$5$-space and a $6$-space are incident if and only if their
intersection is $4$-dimensional.  A $7$-space and a $6$-space are
incident if and only if their intersection is $6$-dimensional (in
which case it is necessarily a nonmaximal singular ideal of $\B$).
\end{thm}

\begin{pf}
Define $X_i$ as in the proof of \ref{e6flag}.  
By Proposition \ref{e7prop}, there is a natural action of
$\Inv(\B)(F_s)$ on 
$X_i(F_s)$, and it is transitive for all $i$ by Theorems \ref{transthm} and
\ref{bigidealthm} and the description of the flag varieties for
$\Inv(J)$ in \ref{e6flag}.

We now show that the roots are associated to the $i$-spaces as
claimed.  We extend scalars so that $\B$ is split.  Let $S$ be a rank
7 split torus in $\Inv(\B)$ as in the proof of Theorem \ref{FTSthm},
where $S_6$ is the torus from the proof of Theorem \ref{e6flag}.

We may extend the characters $\chi_{i,j}$ and $\rho_i$ to $S$ by
setting them to be trivial on $S_1$.  We get a new character $\tau$
defined by
\[
\tau\vert_{S_6} = 1 \hbox{\ and\ } \tau(f_{L_x}) = x.
\]
Then a set of simple roots for $\Inv(\B)$ with respect to $S$ is given
by
\[
\hbox{$\beta_1 := 2 \rho_2 + 2 \tau$ and $\beta_j := \alpha_{j-1}$
for $2 \le j \le 7$.}
\]
For $1 \le j \le 6$, a $j$-space is given by
\[
W_j := \tbtmat{0}{V_{j-1}}{0}{F}
\]
for $V_i$ as in the proof of \ref{e6flag} and $V_0 := 0$.  We also
have a $7$-space
\[
W_7 := \tbtmat{0}{V_6}{Fe_0}{F}
\]
where $e_i$ is the idempotent of $J$ whose only nonzero entry is the
$(i + 1,i + 1)$-entry, which is 1.

The root group for $\beta_1$ is generated by $\psi_{e_2}$ for $\psi$
as in (\ref{psidef}), and the root group for $-\beta_1$ is generated
by $\varphi_{e_2}$ for $\varphi$ as in (\ref{phidef}).

Since $\Inv(\B)(F_s)$ acts transitively on $X_i(F_s)$ and $X_i$ is
clearly a
projective variety for $i \ne 5$, the stabilizer of $W_i$ in
$\Inv(\B)$ is a parabolic subgroup.  For $X_5$, the stabilizer of
$W_5$ is a closed subgroup of $\Inv(\B)$ which contains the Borel
subgroup $P(\emptyset)$ determined by $S$ and our choice of a set of
simple roots.  Thus it is a parabolic subgroup \cite[11.2]{Borel}.
Given the result for $\Inv(J)$, it is now an easy check to see that
the stabilizer of $W_i$ in $\Inv(\B)$ is precisely $P(i)$.  Thus
$X_i = X(i)$.  It only remains to confirm the incidence relations,
which follow easily from the corresponding claims for $E_6$, making
use of the strong restrictions on the form of an inner ideal
containing a singular element observed in the proof of Theorem
\ref{transthm}. 
\end{pf}

It can be hard to visualize what the singular ideals in $\B(J, \D)$
look like.  We mention the particular example  
\begin{neqn} \label{I6}
I_6 := \stbtmat{0}{W}{W}{0} \otimes_F \Delta
\end{neqn}
for $W$ as in (\ref{tredim}).  Then by Lemma \ref{sreglem}, $I_6$ is a 
6-dimensional singular ideal in $\B(J^d, \D)$.  Since it is stable
under $\varpi \otimes \iota$, it is even defined over $F$.

\section{A mysterious result made less so} \label{cisec}
\newcommand{\Esd}{E_6^\D}

There is a nice, but also technical and mysterious, result in
\cite[p.~65, Lem.~3]{Ferr:E6} which describes the form of 1-cocycles
in $Z^1(\D/F, \Esd)$ for $\D$ a quadratic field extension of $F$ and
$\Esd = \autp(\B(J^d, \D))$ the quasi-split simply connected group of
type $E_6$ over $F$ with inner extension $\D$.  Combining Ferrar's
lemma with some newer theorems about groups of type $D_4$ provides
the following stronger result:

\begin{thm} 
For each $\gamma \in H^1(\D/ F, \Esd)$, there is some subgroup $H$ of
$\Esd$ which is simply connected isotropic of type $\oD$ such that
$\gamma$ is in the image of the map $H^1(\D/F, H) \longto H^1(\D/F,
\Esd)$.
\end{thm}

Our proof will not use Ferrar's result nor all of the special Jordan
algebra computations used in its proof.  Our version of the result 
is of additional interest because it can be used to prove that the
Rost invariant $H^1(F, \Esd) \longto H^3(F, \QZt)$ has trivial kernel,
see \cite[\S 31.B]{KMRT} for definitions and a forthcoming paper for a
proof.  We include the theorem here because its proof makes
use of the notations and material in this paper.

\medskip
\begin{pf}
Fix a parabolic subgroup $P = P(\alpha_1, \alpha_6)$ as defined in
\ref{flagback}.  
The action of the nontrivial $F$-automorphism $\iota$
of $\D$ on $\Esd$ is (when the group is considered as a twist of
$\Inv(J^d)$) $^\iota g = \iota g^\dag \iota$, where juxtaposition
denotes the usual $\iota$-action, so $P$ may not be (in fact, isn't)
defined over $F$.
However, the subgroup $G = P \cap {^\iota P}$
certainly is.  Moreover, $\Rel\Cn$ is ``the'' reductive part (= Levi
subgroup) of $P$ and it is defined over $F$, so it is contained in
$G$.  We know $\Rel\Cn = P \cap \op{P}$, where $\op{P}$ is the
opposite parabolic subgroup generated by 
the maximal
torus $S_6$ and the groups $U_\alpha$ for $\alpha \in \pm \D$ such that
$\alpha \ne \alpha_1, \alpha_6$.  We make the following
\begin{neqn} \label{opclaim}
\text{Claim: $^\iota P = \op{P}$,}
\end{neqn}
i.e., that $G = \Rel\Cn$.

We suppose for the moment that the claim is true and return to prove
it later.  A direct consequence of this claim is that the natural map 
\begin{neqn} \label{cohosurj}
H^1(\D/F, G) \longto H^1(\D/F, \Esd)
\end{neqn}
is surjective, as can be seen by making the obvious changes to the
proof of \cite[p.~369, Lem.~6.28]{PlatRap}.

Although $\Rel\Cn$ is defined over $F$, the twisted $\iota$-action on
it is nontrivial.  Since $g_\ut^\dag = g_{\s(\ut)^{-1}}$ for $g$ as in
(\ref{gdef}), the action is given by
\begin{neqn} \label{Reltwist}
^\iota (t_0, t_1, t_2) = \iota (\s(t_0)^{-1}, \s(t_1)^{-1},
\s(t_2)^{-1}) \iota,
\end{neqn}
which restricts to be the usual action of $\iota$ on the subgroup
$\Spin\Cn$.  So what Ferrar proved via Jordan algebra computations in
\cite[p.~65, Lem.~3]{Ferr:E6} was that the map (\ref{cohosurj}) was
surjective.

Suppose now that $\beta \in H^1(\D/F, G)$ is an inverse image of
$\gamma$.  We have an exact sequence over $F_s$
\[
1 \longto \Spin\Cn \longto G \longto K \longto 1
\]
where $K$ is the kernel of the multiplication map
$\mathbb{G}_{m,F}^{\times 3} \longto \GmF$ and the map $G \longto K$
is given by $\trion{t} \mapsto (\mu(t_0), \mu(t_1), \mu(t_2))$.  By
\cite[4.7]{G:iso}, there is a 1-cocycle $a \in Z^1(\D/F, G)$ such that
$(a)$ and $\beta$ have the same image in $H^1(F, K)$ and the twisted
group $\Spin\Cn_a$ is isotropic of type $^1D_4$.  That is, as in
\cite[I.5.5]{SeCG}, we can twist $\Spin\Cn$ by $a$ to obtain a map
\[
\begin{CD}
H^1(\D/F, \Spin\Cn_a) @>>> H^1(\D/F, G_a) @>{\sim}>{\tau_a}>
H^1(\D/F, G)
\end{CD}
\]
which has $\beta$ in its image.  This is the desired result.

We are left with proving Claim (\ref{opclaim}).  The parabolic $\op{P}$ is
obtained from $P$ via the element of the Weyl group denoted by $w_0$
in \cite[p.~220]{Bou:g4}, which is the composition of the map
$\alpha_i \mapsto -\alpha_i$ and the outer automorphism on the Dynkin
diagram.  This shows that $\dim P = \dim \op{P}$.
Since clearly $\dim {^\iota P} = \dim P$, we need only prove that $\op{P}
\subseteq {^\iota P}$.  Since $\Rel\Cn$ (with the twisted
$\iota$-action from (\ref{Reltwist})) is contained in $\op{P}$ and
$^\iota P$, we just need to show that $U_{-\alpha_1}, U_{-\alpha_6} \in
{^\iota P}$.

To decipher the $\iota$-action, we look at the associated root Lie
algebras, which are $\g_{-\alpha_1} = S_{(Fu_1)_{32}}$ and
$\g_{-\alpha_6} = S_{(Fu_1)_{21}}$ in the notation of
\cite[p.~35]{Jac:ex}.  The map $g \mapsto g^\dag$ on $\Inv(J^d)$ has
differential $x \mapsto -x^\ast$.  In general, $S_{u_{ij}}^\ast =
S_{u_{ji}}$, so $\g_{-\alpha_1}^\ast$ is $\g_\gamma$ for 
\[
\gamma = \omega_1 - (\rho_2 - \rho_1) = \alpha_1 + 2\alpha_2 +
2\alpha_3 + \alpha_4 + \alpha_5.
\]
Now $U_\gamma$ is clearly in $P$, so $^\iota U_\gamma = \iota
U_\gamma^\dag \iota = U_{-\alpha_1}$ is in ${^\iota P}$.  A similar
argument works to show that $U_{-\alpha_6}$ is in $^\iota P$, and we
have proven the claim and the theorem.
\end{pf}

\section*{Acknowledgements}
I would like to thank Markus Rost for answering my many questions as well as
for catching an error in an earlier version of this paper.  
Thanks are also due to thank Kevin McCrimmon, Bruce Allison, John Faulkner,
and Oleg Smirnov for providing some references to material on Albert
algebras and 
Ph.~Gille for referring me to the results in Chapter 6 of \cite{PlatRap}.

\providecommand{\bysame}{\leavevmode\hbox to3em{\hrulefill}\thinspace}

\bigskip
%
%
\noindent R.~Skip Garibaldi\\%
e-mail: {\tt skip@member.ams.org} \\%
web: {\tt http://www.math.ucla.edu/\~{}skip/}\\%

\noindent UC Los Angeles\\%
Dept.~of Mathematics\\%
Los Angeles, CA 90095-1555

\end{document}